\documentclass[12pt]{article}

\usepackage
[colorlinks=false,urlcolor=blue,citecolor=blue,linkcolor=blue,bookmarks=true,
     bookmarksopen=false,pdftitle=created by dvipdf,pdfcreator=NM,
     pdfauthor=Megiddo,pdfsubject=mor.sty:6.29.04]{hyperref}

\usepackage{color}
\usepackage{amsfonts}
\usepackage{amssymb}
\usepackage{graphicx}
\usepackage{epsfig}

\def\ignore#1{}

\def\inter{{\rm int}}

\def\ol{\bar}

\def\b{\beta}
\def\d{\delta}
\def\e{\epsilon}

\def\a{\alpha}

\def\tl{\tilde}

\def\rn{\mathbb{R}^n}

\def\re{\mathbb{R}}

\def\dist{{\rm dist}}

\newtheorem{assumption}{Assumption}

\newtheorem{lemma}{Lemma}

\newtheorem{proposition}{Proposition}

\newenvironment{proof}[1][Proof]{\noindent\textbf{#1.} }{\ \rule{0.5em}{0.5em}}
\begin{document}
\title{Constrained Consensus and Optimization in Multi-Agent Networks\footnote{We would like to thank
the associate editor, three anonymous referees, and various seminar
participants for useful comments and discussions.}\ \footnote{This
research was partially supported by the National Science Foundation
under CAREER grants CMMI 07-42538 and DMI-0545910, under grant
ECCS-$0621922$, and by the DARPA ITMANET program}}

\author{Angelia Nedi\'c\footnote{A.\ Nedi\'c is with the
Industrial and Enterprise Systems Engineering Department, University
of Illinois at Urbana-Champaign, Urbana IL 61801
(e-mail: \texttt{angelia@uiuc.edu})}, Asuman Ozdaglar, and Pablo A.\
Parrilo\footnote{A.\ Ozdaglar and P.\ A.\ Parrilo are with the
Laboratory for Information and Decision Systems, Electrical
Engineering and Computer Science Department, Massachusetts Institute
of Technology, Cambridge MA, 02139 (e-mails: \texttt{asuman@mit.edu},
\texttt{parrilo@mit.edu})}}

\markright{LIDS Report 2779}

\maketitle

\thispagestyle{headings}

\begin{abstract}
We present distributed algorithms that can be used by multiple
agents to align their estimates with a particular value over a
network with time-varying connectivity. Our framework is general in
that this value can represent a consensus value among multiple
agents or an optimal solution of an optimization problem, where the
global objective function is a combination of local agent objective
functions. Our main focus is on constrained problems where the
estimate of each agent is restricted to lie in a different
constraint set.

To highlight the effects of constraints, we first consider a
constrained consensus problem and present a distributed ``projected
consensus algorithm" in which agents combine their local averaging
operation with projection on their individual constraint sets. This
algorithm can be viewed as a version of an alternating projection
method with weights that are varying over time and across agents. We
establish convergence and convergence rate results for the projected
consensus algorithm. We next study a constrained optimization
problem for optimizing the sum of local objective functions of the
agents subject to the intersection of their local constraint sets.
We present a distributed ``projected subgradient algorithm" which
involves each agent performing a local averaging operation, taking a
subgradient step to minimize its own objective function, and
projecting on its constraint set. We show that, with an
appropriately selected stepsize rule, the agent estimates generated
by this algorithm converge to the same optimal solution for the
cases when the weights are constant and equal, and when the weights
are time-varying but all agents have the same constraint set.
\end{abstract}

\newpage

\section{Introduction}

There has been much interest in distributed cooperative control
problems, in which several autonomous agents collectively try to
achieve a global objective. Most focus has been on the canonical
consensus problem, where the goal is to develop distributed
algorithms that can be used by a group of agents to reach a common
decision or agreement (on a scalar or vector value). Recent work
also studied multi-agent optimization problems over networks with
time-varying connectivity, where the objective function information
is distributed across agents (e.g., the global objective function is
the sum of local objective functions of agents). Despite much work
in this area, the existing literature does not consider problems
where the agent values are constrained to given sets. Such
constraints are significant in a number of applications including
motion planning and alignment problems, where each agent's position
is limited to a certain region or range, and distributed constrained
multi-agent optimization problems.

In this paper, we study cooperative control problems where the
values of agents are constrained to lie in closed convex sets. Our
main focus is on developing distributed algorithms for problems
where the constraint information is distributed across agents, i.e.,
each agent only knows its own constraint set. To highlight the
effects of different local constraints, we first consider a
constrained consensus problem and propose a {\it projected consensus
algorithm} that operates on the basis of local information. More
specifically, each agent linearly combines its value with those
values received from the time-varying neighboring agents and
projects the combination on its own constraint set.  We show that
this update rule can be viewed as a version of the alternating
projection method where, at each iteration, the values are combined
using weights that are varying in time and across agents, and
projected on the respective constraint sets.

We provide convergence and convergence rate analysis for the
projected consensus algorithm. Due to the projection operation, the
resulting evolution of agent values has nonlinear dynamics, which
poses challenges for the analysis of the algorithm's convergence
properties. To deal with the nonlinear dynamics in the evolution of
the agent estimates, we decompose the dynamics into two parts: a
linear part involving a time-varying averaging operation and a
nonlinear part involving the error due to the projection operation.
This decomposition allows us to represent the evolution of the
estimates using linear dynamics and decouples the analysis of the
effects of constraints from the convergence analysis of the local
agent averaging.  The linear dynamics is analyzed similarly to that
of the unconstrained consensus update, which relies on convergence
of {\it transition matrices} defined as the products of the
time-varying weight matrices.  Using the properties of projection
and agent weights, we prove that the projection error diminishes to
zero. This shows that the nonlinear parts in the dynamics are
vanishing with time and, therefore, the evolution of agent estimates
is ``almost linear''. We then show that the agents reach consensus
on a ``common estimate'' in the limit and that the common estimate
lies in the intersection of the agent individual constraint sets.

We next consider a constrained optimization problem for optimizing a
global objective function which is the sum of local agent objective
functions, subject to a constraint set given by the intersection of
the local agent constraint sets. We focus on distributed algorithms
in which agent values are updated based on local information given
by the agent's objective function and constraint set. In particular,
we propose a distributed {\it projected subgradient algorithm},
which for each agent involves a local averaging operation, a step
along the subgradient of the local objective function, and a
projection on the local constraint set.

We study the convergence behavior of this algorithm for two cases:
when the constraint sets are the same, but the agent connectivity is
time-varying; and when the constraint sets $X_i$ are different, but
the agents use uniform and constant weights in each step, i.e., the
communication graph is fully connected. We show that with an
appropriately selected stepsize rule, the agent estimates generated
by this algorithm converge to the same optimal solution of the
constrained optimization problem. Similar to the analysis of the
projected consensus algorithm, our convergence analysis relies on
showing that the projection errors converge to zero, thus
effectively reducing the problem into an unconstrained one. However,
in this case, establishing the convergence of the projection error
to zero requires understanding the effects of the subgradient steps,
which complicates the analysis. In particular, for the case with
different constraint sets but uniform weights, the analysis uses an
error bound which relates the distances of the iterates to
individual constraint sets with the distances of the iterates to the
intersection set.

Related literature on parallel and distributed computation is vast.
Most literature builds on the seminal work of Tsitsiklis
\cite{johnthes} and Tsitsiklis {\it et al.} \cite{distasyn} (see
also \cite{distbook}), which focused on distributing the
computations involved with optimizing a global objective function
among different processors (assuming complete information about the
global objective function at each processor). More recent literature
focused on multi-agent environments and studied consensus algorithms
for achieving cooperative behavior in a distributed manner (see
\cite{vicsek}, \cite{ali}, \cite{boyd}, \cite{reza},
\cite{spielman}, and \cite{alexCDC, alexlong}). These works assume
that the agent values can be processed arbitrarily and are
unconstrained. Another recent approach for distributed cooperative
control problems involve using game-theoretic models. In this
approach, the agents are endowed with local utility functions that
lead to a game form with a Nash equilibrium which is the same as or
close to a global optimum. Various learning algorithms can then be
used as distributed control schemes that will reach the equilibrium.
In a recent paper, Marden {\it et al.} \cite{JasonECC} used this
approach for the consensus problem where agents have constraints on
their values. Our projected consensus algorithm provides an
alternative approach for this problem.

Most closely related to our work are the recent papers
\cite{ratejournal, ratesubgrad}, which proposed distributed
subgradient methods for solving unconstrained multi-agent
optimization problems. These methods use consensus algorithms as a
mechanism for distributing computations among the agents. The
presence of different local constraints significantly changes the
operation and the analysis of the algorithms, which is our main
focus in this paper. Our work is also related to incremental
subgradient algorithms implemented over a network, where agents
sequentially update an iterate sequence in a cyclic or a random
order \cite{Blatt07,Nedic01,Sundhar08a,Johansson07}. In an
incremental algorithm, there is a single iterate sequence and only
one agent updates the iterate at a given time. Thus, while operating
on the basis of local information, incremental algorithms differ
fundamentally from the algorithm studied in this paper (where all
agents update simultaneously). Furthermore, the work in
\cite{Blatt07,Nedic01,Sundhar08a,Johansson07} assumes that the
constraint set is known by all agents in the system, which is in a
sharp contrast with the algorithms studied in this paper (our
primary interest is in the case where the information about the
constraint set is distributed across the agents).

The paper is organized as follows. In Section~\ref{basics}, we
introduce our notation and terminology, and establish some basic
results related to projection on closed convex sets that will be
used in the subsequent analysis. In Section \ref{sec:consen}, we
present the constrained consensus problem and the projected
consensus algorithm. We describe our multi-agent model and provide a
basic result on the convergence behavior of the transition matrices
that govern the evolution of agent estimates generated by the
algorithms. We study the convergence of the agent estimates and
establish convergence rate results for constant uniform weights.
Section \ref{sec:optim} introduces the constrained multi-agent
optimization problem and presents the projected subgradient
algorithm. We provide convergence analysis for the estimates
generated by this algorithm. Section \ref{conclusions} contains
concluding remarks and some future directions.

\section{Notation, Terminology, and Basics}\label{basics}

\vskip 1pc

A vector is viewed as a column, unless clearly stated otherwise. We
denote by $x_i$ or $[x]_i$ the $i$-th component of a vector $x$.
When $x_i\ge 0$ for all components $i$ of a vector $x$, we write
$x\ge 0$.
%We denote the nonnegative orthant by $\mathbb{R}^m_+$, i.e.,
%$\mathbb{R}^m_+ = \{x\in \mathbb{R}^m\mid x\ge 0\}$.
We write $x'$
to denote the transpose of a vector $x$. The scalar product of two
vectors $x$ and $y$ %\in\re^m$
is denoted by $x'y$. We use $\|x\|$ to denote
the standard Euclidean norm, $\|x\|=\sqrt{x'x}$.
%We write $e_j$ to denote a unit vector with $j$-th component equal to 1
%and all other components equal to 0. We denote by $e$ a vector
%with all components equal to 1.

A vector $a\in\re^m$ is said to be a {\it stochastic vector} when
its components $a_j$ are nonnegative and their sum
is equal to 1, i.e., $\sum_{j=1}^m a_j =1$.
A set of $m$ vectors $\{a^1,\ldots, a^m\}$, %_{i\in \{1,\ldots,m\}}$,
with $a^i\in\re^m$ for all $i$,
is said to be {\it doubly stochastic} when each $a^i$ is
a stochastic vector and $\sum_{i=1}^m a_j^i =1$ for all $j=1,\ldots,m$.
A square matrix
$A$ is said to be doubly stochastic when its rows are stochastic vectors,
and its columns are also stochastic vectors.

We write $dist(\ol x, X)$ to denote the standard
Euclidean distance of a vector $\ol x$ from a set $X$, i.e.,
$$\dist(\ol x,X) =\inf_{x\in X}\|\ol x - x\|.$$
We use $P_X[\ol x]$ to denote the projection of a vector $\ol x$
on a closed convex set $X$, i.e.,
\[P_X[\ol x]=\arg\min_{x\in X}\|\ol x-x\|.\]
In the subsequent development, the properties of the projection
operation on a closed convex set play an important role. In
particular, we use the projection inequality, i.e., for any vector
$x$,
\begin{equation}\label{ooo}
(P_X[x]-x)'\,(y-P_X[x])\ge 0\qquad
\hbox{for all }y\in X.
\end{equation}
We also use
the standard non-expansiveness property, i.e.,
\begin{equation}\label{nonexpan}
\|P_X[x] - P_X[y]\|\le \|x-y\|\qquad\hbox{for any }x \hbox{ and } y.
\end{equation}
In addition, we use the properties given in the following lemma.

\begin{lemma}\label{proj_prop}\em{
Let $X$ be a nonempty closed convex set in $\rn$. Then,
we have for any $x\in\rn$,
\begin{itemize}
\item[(a)]
\qquad$(P_X[x]-x)'\,(x-y)\le -\|P_X[x]-x\|^2\qquad\hbox{for all $y\in X$}$.
\item[(b)]
\qquad$\|P_X[x]-y\|^2\le \|x-y\|^2 - \|P_X[x]-x\|^2
\qquad\hbox{for all $y\in X$}$.
\end{itemize}
}
\end{lemma}

\begin{proof}
(a)\ Let $x\in\rn$ be arbitrary. Then, for any $y\in X$, we have
\[
(P_X[x]-x)'\,(x-y)
= (P_X[x]-x)'\,(x-P_X[x])
+(P_X[x]-x)'\, (P_X[x]- y).\]
By the projection inequality [cf.\ Eq.\ (\ref{ooo})], it follows that
$(P_X[x]-x)'\, (P_X[x]- y)\le 0$,
implying
\[(P_X[x]-x)'\,(x-y)\le -\|P_X[x]-x\|^2
\qquad\hbox{for all }y\in X.\]

\noindent (b)\ For an arbitrary $x\in\rn$ and for all $y\in X$, we
have
\begin{eqnarray*}
\|P_X[x]-y\|^2 &=& \|P_X[x]-x+x-y\|^2\\
&=& \|P_X[x]-x\|^2 + \|x-y\|^2+
2(P_X[x]-x)'(x-y).
\end{eqnarray*}
By using the inequality of part (a),
we obtain
$$\|P_X[x]-y\|^2\le \|x-y\|^2 - \|P_X[x]-x\|^2
\qquad\hbox{for all }y\in X.$$
\end{proof}

\begin{figure}\label{projfigure}
\begin{center}
\epsfxsize=3in \epsffile{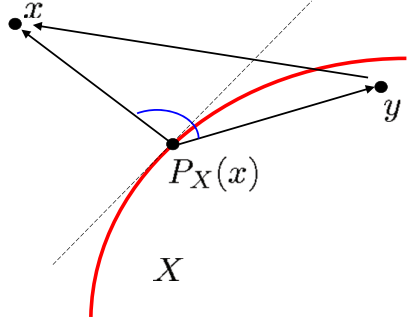}
\caption{\small Illustration of the relation between the projection error
and feasible directions of a convex set at the projection vector.}
\end{center}
\end{figure}

\vskip 1pc Part (b) of the preceding Lemma establishes a relation
between the projection error vector and the feasible directions of
the convex set $X$ at the projection vector, as illustrated in
Figure \ref{projfigure}.

We next consider nonempty closed convex sets $X_i\subseteq \rn$, for
$i=1,\ldots,m$, and an averaged-vector $\hat x$ obtained by taking
an average of vectors $x_i\in X_i$, i.e., $\hat x={1\over m}
\sum_{i=1}^m x_i$ for some $x_i\in X_i$. We provide an ``error
bound" that relates the distance of the  averaged-vector $\hat x$
from the intersection set $X=\cap_{i=1}^m X_i$ to the distance of
$\hat x$ from the individual sets $X_i$. This relation, which is
also of independent interest, will play a key role in our analysis
of the convergence of projection errors associated with various
distributed algorithms introduced in this paper. We establish the
relation under an interior point assumption on the intersection set
$X=\cap_{i=1}^m X_i$ stated in the following:

\begin{assumption}
\emph{Given sets $X_i\subseteq \rn$, $i=1,\ldots,m$, let
$X=\cap_{i=1}^m X_i$ denote their intersection.
There is a vector %$\bar x\in X$ such that
$\bar x \in \inter(X),$ %=\inter\Big(\cap_{i=1}^m X_i\Big),\]
i.e., there exists a scalar $ \delta>0$ such that
$$\{z\ |\ \|z-\bar x\|\le \delta\}\subset X.$$
}\label{assump:intcond}
\end{assumption}

We provide an error bound relation in the following lemma.

\begin{lemma}
\emph{Let $X_i\subseteq \rn$, $i=1,\ldots,m$, be nonempty closed
convex sets that satisfy Assumption \ref{assump:intcond}. Let
$x^i\in X_i$, $i=1,\ldots,m$, be arbitrary vectors and define their
average as $\hat x={1\over m} \sum_{i=1}^m x^i$. Consider the vector
$s\in \rn$ defined by
\[s={\e\over \e+\delta}\, \bar x + {\delta \over \e+\delta}\, \hat x,\]
where
\[\e = \sum_{j=1}^m \dist(\hat x, X_j),\] and $\delta$ is the
scalar given in Assumption \ref{assump:intcond}.
\begin{itemize}
\item[(a)] The vector $s$
belongs to the intersection set $X=\cap_{i=1}^m X_i$.
\item[(b)] We have the following relation
\[\|\hat x - s\| \le {1\over \delta m} \Big(\sum_{j=1}^m \|x^j-\bar
x\|\Big) \Big(\sum_{j=1}^m \dist(\hat x, X_j)\Big).\] As a
particular consequence, we have
\[\dist(\hat x,X) \le {1\over \delta m} \Big(\sum_{j=1}^m \|x^j-\bar x\|\Big) \Big(\sum_{j=1}^m \dist(\hat x, X_j)\Big).\]
\end{itemize}
}\label{distance}
\end{lemma}

\begin{proof}

\noindent (a)\ We first show that the vector $s$ belongs to the
intersection $X=\cap_{i=1}^m X_i$. To see this, let
$i\in\{1,\ldots,m\}$ be arbitrary and note that we can write $s$ as
\[s={\e\over \e+\delta}
\left(\bar x + {\delta\over \e} \Big(\hat x-P_{X_i}[\hat
x]\Big)\right) +{\delta\over \e+\delta}P_{X_i} [\hat x].\] By the
definition of $\e$, it follows that $\|\hat x-P_{X_i}[\hat x]\|\le
\e$, implying by the interior point assumption (cf.\ Assumption
\ref{assump:intcond}) that the vector $\bar x + {\delta\over \e}
\Big(\hat x-P_{X_i}[\hat x]\Big)$ belongs to the set $X$, and
therefore to the set $X_i$. Since the vector $s$ is the convex
combination of two vectors in the set $X_i$, it follows by the
convexity of $X_i$ that $s\in X_i$. The preceding argument is valid
for an arbitrary $i$, thus implying that $s\in X$.

\vskip 1pc

\noindent (b)\ Using the  definition of the vector $s$ and the
vector $\hat x$, we have
\begin{eqnarray*}
\|\hat x -s\| = \frac{\e}{\e+\d} \left\|{1\over m} \sum_{j=1}^m x^j
-\bar x \right\|\le \frac{\e}{\d m}\,\sum_{j=1}^m   \|x^j -\bar x
\|.
\end{eqnarray*}
Substituting the definition of $\e$ yields the desired relation.
\end{proof}

\section{Constrained Consensus}\label{sec:consen}

In this section, we describe the constrained consensus problem. In
particular, we introduce our multi-agent model and the projected
consensus algorithm that is locally executed by each agent. We
provide some insights about the algorithm and we discuss its
connection to the alternating projections method. We also introduce
the assumptions on the multi-agent model and present key elementary
results that we use in our subsequent analysis of the projected
consensus algorithm. In particular, we define the transition
matrices governing the linear dynamics of the agent estimate
evolution and give a basic convergence result for these %transition
matrices. The model assumptions and the transition matrix
convergence properties will also be used for studying the
constrained optimization problem and the projected subgradient
algorithm that we introduce in Section \ref{sec:optim}.

\subsection{Multi-Agent Model and Algorithm}\label{model_algo}

We consider a set of agents denoted by $V=\{1,\ldots,m\}.$ We assume
a slotted-time system, and we denote by $x^i(k)$ the estimate
generated and stored by agent $i$ at time slot $k$. The agent
estimate $x^i(k)$ is a vector in $\rn$ that is constrained to lie in
a nonempty closed convex set $X_i\subseteq\rn$ known only to agent
$i$. The agents' objective is to cooperatively  reach a consensus on
a common vector through a sequence of local estimate updates
(subject to the local constraint set) and local information
exchanges (with neighboring agents only).

We study a model where the agents exchange and update their
estimates as follows: To generate the estimate at time $k+1$, agent
$i$ forms a convex combination of its estimate $x^i(k)$ with the
estimates received from other agents at time $k$, and takes the
projection of this vector on its constraint set $X_i$. More
specifically, agent $i$ at time $k+1$ generates its new estimate
according to the following relation:
\begin{equation}
x^i(k+1) = P_{X_i}\left[\sum_{j=1}^m a_j^i(k)
x^j(k)\right],\label{proj-estimate}
\end{equation}
where $a^i=(a_1^i,\ldots,a_m^i)'$ is a vector of nonnegative
weights.

The relation in Eq.\ (\ref{proj-estimate}) defines the {\it
projected consensus algorithm}. The method can be interpreted as a
multi-agent algorithm for finding a point in common to the given
closed convex sets $X_1,\ldots, X_m$. Note that the problem of
finding a common point can be formulated as an unconstrained convex
optimization problem of the following form:
\begin{equation}\label{commonpoint}
\begin{array} {ll}
\hbox{minimize } & \frac{1}{2}\, \sum_{i=1}^m
\left\|x-P_{X_i}[x]\right\|^2\cr \hbox{subject to } & x\in\rn.
\end{array}
\end{equation}
In view of this optimization problem, the method can be interpreted
as a {\it distributed gradient algorithm} where each agent is
assigned an objective function
$f_i(x)=\frac{1}{2}\,\left\|x-P_{X_i}[x]\right\|^2$. At each time
$k+1$, an agent incorporates new information $x^j(k)$ received from
some of the other agents and generates a weighted sum $\sum_{j=1}^m
a^i_j(k)x^j(k)$. Then, the agent updates its estimate by taking a
step (with stepsize equal to 1) along the negative gradient of its
own objective function $f_i=\frac{1}{2}\,\|x-P_{X_i}\|^2$ at
$x=\sum_{j=1}^m a^i_j(k)x^j(k)$. In particular, since the gradient
of $f_i$ is $\nabla f_i(x) =x-P_{X_i}[x]$ (see Theorem 1.5.5 in
Facchinei and Pang \cite{facc_pang}), the update rule in Eq.\
(\ref{proj-estimate}) is equivalent to the following gradient
descent method for minimizing $f_i$:
\[x^i(k+1) = \sum_{j=1}^m a_j^i(k)x^j(k) -\left(
\sum_{j=1}^m a_j^i(k)x^j(k)  - P_{X_i} \left[\sum_{j=1}^m
a_j^i(k)x^j(k)\right]\right).
\]
This view of the update rule motivates our line of  analysis of the
projected consensus method. In particular, motivated by the
objective function of problem (\ref{commonpoint}), we use
$\sum_{i=1}^m \left\|x^i(k)-x\right\|^2$ with $x\in \cap_{i=1}^m
X_i$ as a Lyapunov function measuring the progress of the algorithm
(see Section \ref{rate_analysis}).\footnote{ We focus throughout the
paper on the case when the intersection set $\cap_{i=1}^m X_i$ is
nonempty. If the intersection set is empty, it follows from the
definition of the algorithm that the agent estimates will not reach
a consensus. In this case, the estimate sequences $\{x^i(k)\}$ may
exhibit oscillatory behavior or may all be unbounded.}

\subsection{Relation to Alternating Projections Method}\label{relaltproj}

\begin{figure}
\centerline{\hbox{ \hspace{0.0in}
    \epsfxsize=3.5in
    \epsffile{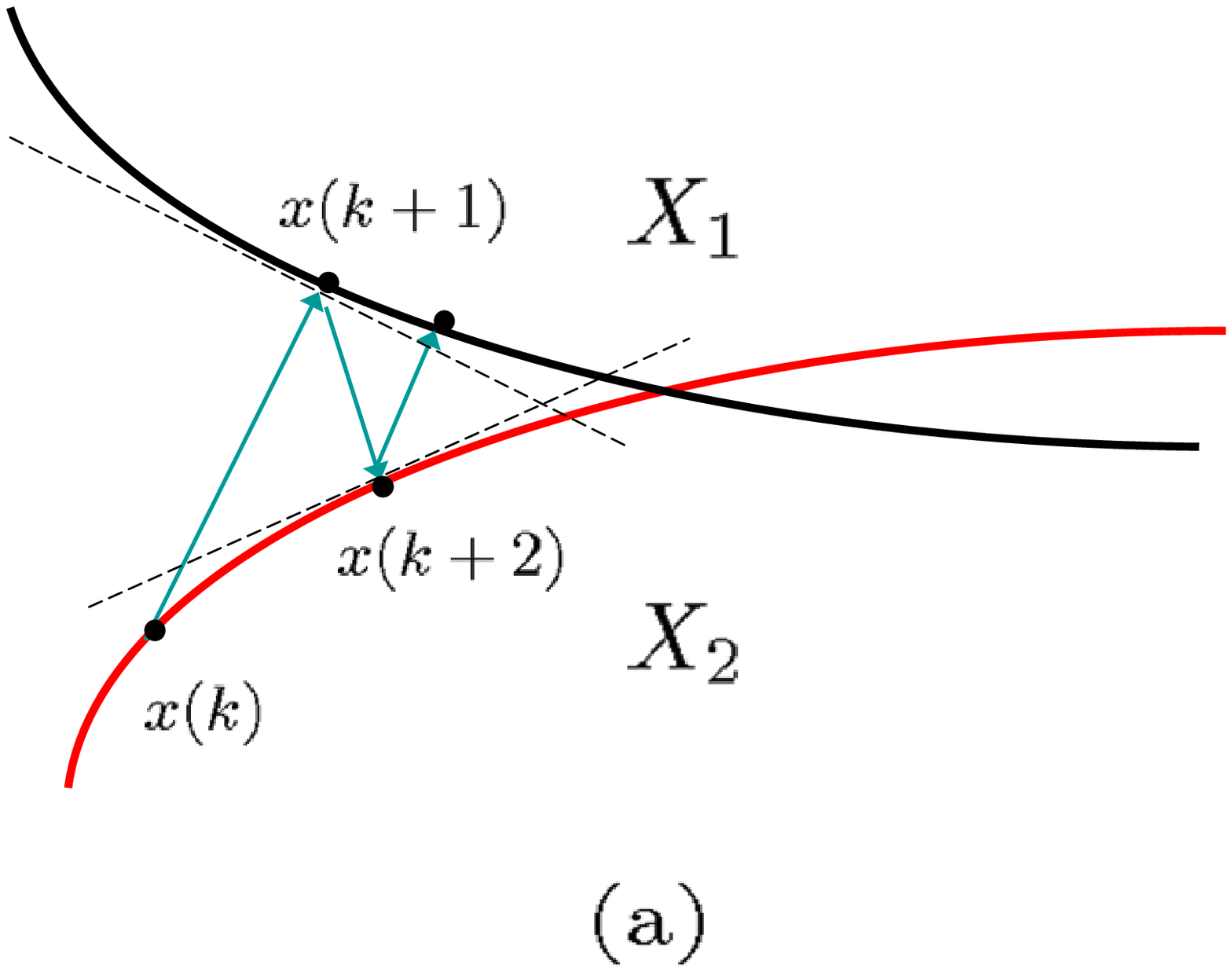}
\hfill\hfill
    \epsfxsize=3.5in
   \epsffile{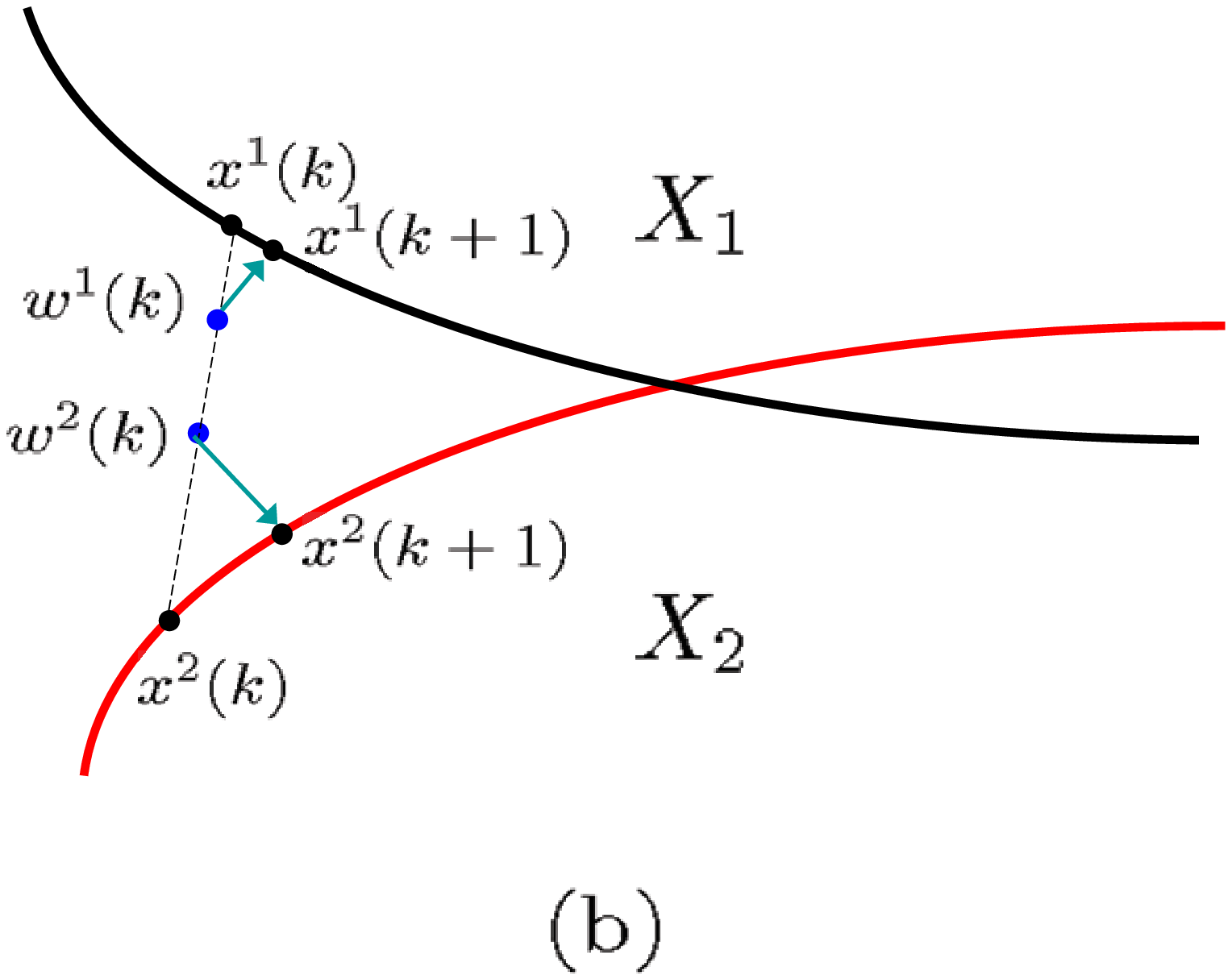}
    }
  }
\caption{\small Illustration of the connection between the
alternating/cyclic projection method and the constrained consensus
algorithm for two closed convex sets $X_1$ and $X_2$. In plot (a),
the alternating projection algorithm generates a sequence $\{x(k)\}$
by iteratively projecting onto sets $X_1$ and $X_2$, i.e.,
$x(k+1)=P_{X_1}[x(k)],\ x(k+2)=P_{X_2}[x(k+1)]$. In plot (b), the
projected consensus algorithm generates sequences $\{x^i(k)\}$ for
agents $i=1,2$ by first combining the iterates with different
weights and then projecting on respective sets $X_i$, i.e.,
$w^i(k)=\sum_{j=1}^m a_j^i(k)x^j(k)$ and $x^i(k+1)=P_{X_i}[w^i(k)]$
for $i=1,2$. } \label{proj-methods}
\end{figure}

The method of Eq.\ (\ref{proj-estimate}) is related to the classical
{\it alternating or cyclic projection method}. Given a finite
collection of closed convex sets $\{X_i\}_{i\in {\cal I}}$ with a
nonempty intersection (i.e., $\cap_{i\in {\cal I}}X_i\ne\emptyset$),
the alternating projection method finds a vector in the intersection
$\cap_{i\in {\cal I}}X_i$. In other words, the algorithm solves the
unconstrained problem (\ref{commonpoint}). Alternating projection
methods generate a sequence of vectors by projecting iteratively on
the sets (either cyclically or with some given order), see Figure
\ref{proj-methods}(a). The convergence behavior of these methods has
been established by Von Neumann \cite{vonNeumann} and Aronszajn
\cite{aronszajn} for the case when the sets $X_i$ are affine; and by
Gubin {\it et al.} \cite{gubinpolyak} when the sets $X_i$ are closed
and convex. Gubin {\it et al.} \cite{gubinpolyak} also have provided
convergence rate results for a particular form of alternating
projection method. Similar rate results under different assumptions
have also been provided by Deutsch \cite{deutsch}, and Deutsch and
Hundal \cite{deuhun}.

The constrained consensus algorithm [cf.\ Eq.\
(\ref{proj-estimate})] generates a sequence of iterates for each
agent as follows: at iteration $k$, each agent $i$ first forms a
linear combination of the other agent values $x^j(k)$ using its own
weight vector $a^i(k)$ and  then projects this combination on its
constraint set $X_i$. Therefore, the projected consensus algorithm
can be viewed as a version of the alternating projection algorithm,
where the iterates are combined with the weights varying over time
and across agents, and then projected on the individual constraint
sets, see Figure \ref{proj-methods}(b).

We conclude this section by noting that the alternate projection
method has much more structured weights than the weights we consider
in this paper. As seen from the assumptions on the agent weights in
the following section, the analysis of our projected consensus
algorithm (and the projected subgradient algorithm introduced in
Section \ref{sec:optim}) is complicated by the general time
variability of the weights $a^i_j(k)$.

\subsection{Assumptions}\label{sec:assumptions}
Following Tsitsiklis \cite{johnthes} (see also Blondel {\it et al.}
\cite{multiagent}), we adopt the following assumptions on the weight
vectors $a^i(k)$, $i\in \{1,\ldots,m\}$ and on information exchange.

\begin{assumption} (Weights Rule)
\emph {There exists a scalar $\eta$ with $0<\eta< 1$ such that
for all $i\in\{1,\ldots,m\}$,
\begin{itemize}
\item[(a)] $a_i^i(k) \ge \eta$ for all $k\ge0$.
\item[(b)] If $a_j^i(k)>0$, then $a_j^i(k)\ge\eta$.
\end{itemize}
}\label{weightsrule}
\end{assumption}

\begin{assumption}(Doubly Stochasticity) \emph{
The vectors $a^i(k)=(a^i_1(k),\ldots,a^i_m(k))'$ satisfy:
\begin{itemize}
\item[(a)] $a^i(k)\ge0$
and $\sum_{j=1}^m a_j^i(k) = 1$ for all $i$ and $k$, i.e.,
the vectors $a^i(k)$ are stochastic.
\item[(b)] $\sum_{i=1}^m a^i_j(k) = 1$ for all $j$ and $k$.
\end{itemize}
} \label{doubstoc}
\end{assumption}

Informally speaking, Assumption \ref{weightsrule} says that every
agent assigns a substantial weight to the information received from
its neighbors. This guarantees that the information from each agent
influences the information of every other agent persistently in
time. In other words, this assumption guarantees that the agent
information is mixing at a nondiminishing rate in time. Without this
assumption, information from some of the agents may become less
influential in time, and in the limit, resulting in loss of
information from these agents.

Assumption \ref{doubstoc}(a) establishes that each agent takes a
convex combination of its estimate and the estimates of its
neighbors. Assumption \ref{doubstoc}(b), together with Assumption
\ref{weightsrule}, ensures that the estimate of every agent is
influenced by the estimates of every other agent with the same
frequency in the limit, i.e., all agents are equally influential in
the long run.

We now impose some rules on the agent information exchange.
At each update time $t_k$,
the information exchange among the agents may be represented by
a directed graph $(V,E_k)$ with the set $E_k$
of directed edges given by
$$E_k=\{(j,i)\mid a^i_j(k)>0\}.$$
Note that, by Assumption \ref{weightsrule}(a),
we have $(i,i)\in E_k$ for each agent $i$ and all $k$. Also, we have
$(j,i)\in E_k$
if and only if agent $i$ receives the information $x^j$ from agent $j$
in the time interval $(t_k,t_{k+1})$.

We next formally state the connectivity assumption on the
multi-agent system. This assumption ensures that the information of
any agent $i$ influences the information state of any other agent
infinitely often in time.

\begin{assumption} \label{connected}
(Connectivity) \emph{ The graph $(V, E_\infty)$ is strongly
connected, where $E_\infty$  is the set of edges $(j,i)$
representing agent pairs communicating directly infinitely many
times, i.e.,
$$E_\infty=\{(j,i)\mid
(j,i)\in E_k \hbox{ for infinitely many indices }k \}.$$ }
\end{assumption}

We also adopt an additional
assumption that
the intercommunication intervals are bounded for those agents that
communicate directly. In particular, this is stated in the following.

\begin{assumption} (Bounded Intercommunication Interval) \emph {There
exists an integer $B\ge 1$ such that for every $(j,i)\in E_\infty$,
agent $j$ sends its information to a neighboring agent $i$ at
least once every $B$ consecutive time slots, i.e., at time $t_k$ or
at time $t_{k+1}$ or $\ldots$ or (at latest) at time $t_{k+B-1}$ for
any $k\ge0$. }\label{boundedintervals}
\end{assumption}

In other words, the preceding assumption guarantees that every pair
of agents that communicate directly infinitely many times exchange
information at least once every $B$ time slots.\footnote{It is
possible to adopt weaker connectivity assumptions for the
multi-agent model as those used in the recent work
\cite{alex_quant_pap}.}

\subsection{Transition Matrices}\label{tranmat}
We introduce matrices $A(s)$, whose $i$-th column is the weight
vector $a^i(s)$, and the matrices
\[
\Phi(k,s) = A(s)A(s+1)\cdots A(k-1)A(k)\qquad \hbox{for all $s$ and
$k$ with $k\ge s$,}
\]
where
\[
\Phi(k,k) = A(k)\qquad \hbox{for all $k$}.
\]
We use these matrices to describe the evolution of the agent
estimates associated with the algorithms introduced in Sections
\ref{sec:consen} and \ref{sec:optim}. The convergence properties of
these matrices as $k\to \infty$ have been extensively studied and
well-established (see \cite{johnthes}, \cite{ali}, \cite{wolf}).
 Under the assumptions of Section
\ref{sec:assumptions}, the matrices $\Phi(k,s)$ converge as
$k\to\infty$ to a {\it uniform steady state distribution for each
$s$ at a geometric rate}, i.e., $\lim_{k\to\infty}\Phi(k,s)= {1\over
m}\, ee'$ for all $s$. The fact that transition matrices converge at
a geometric rate plays a crucial role in our analysis of the
algorithms. Recent work has established explicit convergence rate
results for the transition matrices \cite{ratejournal,ratesubgrad}.
These results are given in the following proposition without a
proof.

\begin{proposition}
\emph{Let Assumptions \ref{weightsrule}, \ref{doubstoc},
\ref{connected} and \ref{boundedintervals} hold. Then, we have the
following:
\begin{itemize}
\item[(a)]
The entries $[\Phi(k,s)]^i_j$ of the transition matrices
converge to ${1\over m}$ as
$k\to\infty$ with a geometric rate uniformly with respect to $i$ and
$j$, i.e., for all $i,j\in\{1,\ldots,m\}$,
$$\left|[\Phi(k,s)]^i_j - {1\over m}\right|
\le 2\,{1+\eta^{-B_0}\over 1-\eta^{B_0}}\
\left(1-\eta^{B_0}\right)^{{k-s\over B_0}} \quad \hbox{for all $s$
and $k$ with }k\ge s.$$
\item[(b)]
In the absence of Assumption \ref{doubstoc}(b) [i.e., the weights
$a^i(k)$ are stochastic but not doubly stochastic], the columns
$[\Phi(k,s)]^i$ of the transition matrices converge to a stochastic
vector $\phi(s)$ as $k\to\infty$ with a geometric rate uniformly
with respect to $i$ and $j$, i.e., for all $i,j\in\{1,\ldots,m\}$,
$$\left|[\Phi(k,s)]^i_j - \phi_j(s)\right|
\le 2\,{1+\eta^{-B_0}\over 1-\eta^{B_0}}\
\left(1-\eta^{B_0}\right)^{{k-s\over B_0}} \quad \hbox{for all $s$
and $k$ with }k\ge s.$$
\end{itemize}
Here, $\eta$ is the lower bound of
Assumption \ref{weightsrule}, $ B_0 = (m-1)B$, $m$ is the number
of agents, and $B$ is the intercommunication interval bound of
Assumption \ref{boundedintervals}.
}\label{tranconv}
\end{proposition}

\subsection{Convergence}\label{convergence} In this section, we study
the convergence behavior of the agent estimates $\{x^i(k)\}$
generated by the projected consensus algorithm (\ref{proj-estimate})
under Assumptions \ref{weightsrule}--\ref{boundedintervals}. We
write the update rule in Eq.\ (\ref{proj-estimate}) as
\begin{equation}
x^i(k+1) =\sum_{j=1}^m a_j^i(k) x^j(k)+e^i(k),
\label{decomp_dyn}
\end{equation}
where $e^i(k)$ represents the error due to projection given by
\begin{equation}
e^i(k)= P_{X_i}\left[\sum_{j=1}^m a_j^i(k) x^j(k)\right]
-\sum_{j=1}^m a_j^i(k) x^j(k).
\label{nonlin_error}
\end{equation}

As indicated by the preceding two relations, the evolution dynamics
of the estimates $x^i(k)$ for each agent is decomposed into a sum of
a linear (time-varying) term $\sum_{j=1}^m a_j^i(k) x^j(k)$ and a
nonlinear term $e^i(k)$. The linear term captures the effects of
mixing the agent estimates, while the nonlinear term captures the
nonlinear effects of the projection operation. This decomposition
plays a crucial role in our analysis. As we will shortly see [cf.\
Lemma \ref{seqbehav}(d)], under the doubly stochasticity assumption
on the weights, the nonlinear terms $e^i(k)$ are diminishing in time
for each $i$, and therefore, the evolution of agent estimates is
``almost linear''. Thus, the nonlinear term can be viewed as a
non-persistent disturbance in the linear evolution of the estimates.

For notational convenience, let $w^i(k)$ denote
\begin{equation}
w^i(k)= \sum_{j=1}^m a_j^i(k) x^j(k).\label{wvector}\end{equation}
In this notation, the iterate $x^i(k+1)$ and the projection error
$e^i(k)$ are given by
\begin{equation}
x^i(k+1)=P_{X_i} [w^i(k)], \label{error-estimate}
\end{equation}
\begin{equation}
e^i(k)=x^i(k+1)-w^i(k).\label{error}
\end{equation}

In the following lemma, we show some relations for the sums
$\sum_{i=1}^m\|x^i(k)-x\|^2$ and $\sum_{i=1}^m\|w^i(k)-x\|^2$, and
$\sum_{i=1}^m\|x^i(k)-x\|$ and $\sum_{i=1}^m\|w^i(k)-x\|$ for an
arbitrary vector $x$ in the intersection of the agent constraint
sets. Also, we  prove that the errors $e^i(k)$
converge to zero as $k\to\infty$ for all $i$.
The projection properties given in Lemma \ref{proj_prop} and
the doubly stochasticity of the weights play crucial roles in
establishing these relations.

\begin{lemma}
\emph{Let the intersection set $X=\cap_{i=1}^m X_i$ be nonempty, and
let Doubly Stochasticity assumption hold (cf.\ Assumption
\ref{doubstoc}). Let $x^i(k)$, $w^i(k)$, and $e^i(k)$ be defined by
Eqs.\ (\ref{wvector})--(\ref{error}). Then, we have the following.
\begin{itemize}
\item[(a)]
For all $x\in X$ and all $k$, we have
\begin{itemize}
\item[(i)]\quad
$\|x^i(k+1)-x\|^2\le\|w^i(k)-x\|^2-\|e^i(k)\|^2\qquad\hbox{for all }i,$
\item[(ii)]
\quad $\sum_{i=1}^m \|w^i(k)-x\|^2 \le
\sum_{i=1}^m \|x^i(k)-x\|^2,$
\item[(iii)]\quad
$\sum_{i=1}^m \|w^i(k)-x\| \le \sum_{i=1}^m \|x^i(k) - x\|$.
\end{itemize}
\item[(b)]
For all $x\in X$, the sequences $\Big\{\sum_{i=1}^m
\|w^i(k)-x\|^2\Big\}$ and $\Big\{\sum_{i=1}^m \|x^i(k)-x\|^2\Big\}$
are monotonically nonincreasing with $k$.
\item[(c)]
For all $x\in X$, the sequences $\Big\{\sum_{i=1}^m
\|w^i(k)-x\|\Big\}$ and $\Big\{\sum_{i=1}^m \|x^i(k)-x\|\Big\}$ are
monotonically nonincreasing with $k$.
\item[(d)]
The errors $e^i(k)$ converge to zero as $k\to\infty$, i.e.,
$$\lim_{k\to\infty} e^i(k)=0\qquad\hbox{for all }i.$$
\end{itemize}
}\label{seqbehav}\end{lemma}

\begin{proof}
(a)\ For any $x\in X$ and $i$, we consider the term $\|x^i(k+1)
-x\|^2$. Since $X\subseteq X_i$ for all $i$, it follows that $x\in
X_i$ for all $i$. Since we also have $x^i(k+1)=P_{X_i} [w^i(k)]$, we
have from Lemma \ref{proj_prop}(b) that
\[
\left\|x^i(k+1) -x\right\|^2 \le \left\|w^i(k) -x \right\|^2
-\|x^i(k+1)-w^i(k)\|^2\qquad \hbox{for all }x\in X\hbox{ and }k\ge 0,
\]
which yields the relation in part (a)(i) in view of relation
(\ref{error}).

By the definition of $w^i(k)$ in Eq.\ (\ref{wvector}) and the
stochasticity of the weight vector $a^i(k)$ [cf.\ Assumption
\ref{doubstoc}(a)], we have for every agent $i$ and any $x\in X$,
\begin{equation}
w^i(k)-x=\sum_{j=1}^m a^i_j(k)\left(x^j(k)-x\right) \qquad\hbox{for
all }k\ge0.\label{convcomb}
\end{equation}
Thus, for any $x\in X$, and all $i$ and $k$,
\[
\|w^i(k)-x\|^2
=\left\|\sum_{j=1}^m a^i_j(k)\left(x^j(k)-x\right)\right\|^2
\le\sum_{j=1}^m a^i_j(k)\left\|x^j(k) -x\right\|^2,
\]
where the inequality holds since
the vector $\sum_{j=1}^m
a^i_j(k)(x^j(k)-x)$ is a convex combination of the vectors
$x^j(k)-x$ and the squared norm
$\|\cdot\|^2$ is a convex function.
By summing the preceding relations over $i=1,\ldots,m$,
we obtain
\[
\sum_{i=1}^m \|w^i(k)-x\|^2  \le \sum_{i=1}^m \sum_{j=1}^m
a^i_j(k)\left\|x^j(k) -x\right\|^2
= \sum_{j=1}^m
\left(\sum_{i=1}^m a^i_j(k)\right) \left\|x^j(k) -x\right\|^2.
\]
Using the doubly stochasticity of the weight vectors $a^i(k)$, i.e.,
$\sum_{i=1}^m a_j^i(k)=1$ for all $j$ and $k$
[cf.\ Assumption \ref{doubstoc}(b)], we obtain the relation in part (a)(ii),
\[%\begin{equation}
\sum_{i=1}^m \|w^i(k)-x\|^2 \le \sum_{i=1}^m \left\|x^i(k)
-x\right\|^2\qquad \hbox{for all }x\in X
\hbox{ and }k\ge0.%\label{dstocdec}
\]%\end{equation}

Similarly, from relation (\ref{convcomb}) and the doubly
stochasticity of the weights, we obtain for all $x\in X$ and all
$k$,
\[
\sum_{i=1}^m \|w^i(k)-x\| \le \sum_{i=1}^m \sum_{j=1}^m a_j^i(k)
\|x^j(k) - x\| = \sum_{j=1}^m \|x^j(k) - x\|, %\label{monds}
\]
thus showing the relation in part (a)(iii).

\vskip 0.5pc

\noindent (b)\ For any $x\in X$, the nonincreasing properties of the
sequences $\Big\{\sum_{i=1}^m \|w^i(k)-x\|^2\Big\}$ and
$\Big\{\sum_{i=1}^m \|x^i(k)-x\|^2\Big\}$ follow by
combining the relations in parts (a)(i)--(ii).
%(\ref{xiwiei}) and (\ref{dstocdec}).

\vskip 0.5pc
\noindent (c)\ Since $x^i(k+1)=P_{X_i}(w^i(k))$ for all $i$ and
$k\ge 0$, using the nonexpansiveness property of the projection
operation [cf.\ Eq.\ (\ref{nonexpan})], we have
\[
\|x^i(k+1)-x\|\le \|w^i(k)-x\|\qquad \hbox{for all }x\in X_i,\
\hbox{ all } i \hbox{ and } k.
\]
Summing the preceding relations over all $i\in \{1,\ldots,m\}$ yields
for all $k$,
\begin{equation}
\sum_{i=1}^m \|x^i(k+1)-x\|\le \sum_{i=1}^m \|w^i(k)-x\|\qquad
\hbox{for all }x\in X.\label{nonexp}
\end{equation}
The nonincreasing property of the sequences
$\{\sum_{i=1}^m \|w^i(k)-x\|\}$ and $\{\sum_{i=1}^m
\|x^i(k)-x\|\}$ follows
from the preceding relation and the relation in part (a)(iii).
% (\ref{nonexp}) and (\ref{monds}).
\vskip 0.5pc
\noindent
(d)\ By summing the relations
in part (a)(i)
over $i=1,\ldots,m$, we obtain for any
$x\in X$,
\[\sum_{i=1}^m \left\|x^i(k+1) -x\right\|^2
\le \sum_{i=1}^m \left\|w^i(k) -x \right\|^2
-\sum_{i=1}^m \|e^i(k)\|^2\qquad \hbox{for all }k\ge 0.\]
Combined
with the inequality $\sum_{j=1}^m \|w^j(k)-x\|^2\le
\sum_{j=1}^m \|x^j(k)-x\|^2$ of part (a)(ii),
we further obtain
\[%\begin{equation}
\sum_{i=1}^m \|e^i(k)\|^2\le
\sum_{i=1}^m \left\|x^i(k) -x \right\|^2 -
\sum_{i=1}^m \left\|x^i(k+1) -x\right\|^2 \qquad \hbox{for all }k\ge 0.
%\label{exkminus}
\]%\end{equation}
Summing these relations over $k=0,\ldots, s$ for any
$s>0$ yields
\[\sum_{k=0}^s\,\sum_{i=1}^m\|e^i(s)\|^2
\le \sum_{i=1}^m\,\left\|x^i(0) -x\right\|^2
-\sum_{i=1}^m \left\|x^i(s+1) -x\right\|^2
\le \sum_{i=1}^m\,\left\|x^i(0) -x\right\|^2.
%\qquad\hbox{for any }x\in X.
\]
By  letting $s\to\infty$, we obtain
$$\sum_{k=0}^\infty\,\sum_{i=1}^m\|e^i(k)\|^2
\le\sum_{i=1}^m\,\left\|x^i(0) -x\right\|^2,$$
implying
$\lim_{k\to\infty}\|e^i(k)\|=0$ for all $i$.
\end{proof}

\vskip 1pc

We next consider the evolution of the estimates
$x^i(k+1)$ generated by method (\ref{proj-estimate})
over a period of time.
In particular, we relate the estimates $x^i(k+1)$ to the
estimates $x^i(s)$ generated earlier in time $s$ with $s<k+1$
by exploiting the decomposition of
the estimate evolution in Eqs.\ (\ref{decomp_dyn})--(\ref{nonlin_error}).
In this, we use the transition matrices $\Phi(k,s)$ from
time $s$ to time $k$ (see Section \ref{tranmat}). As we will shortly see,
the linear part of the dynamics is given in
terms of the transition matrices, while the
nonlinear part involves combinations
of the transition matrices and the error terms from
time $s$ to time $k$.

Recall that the transition matrices are defined as follows:
\[
\Phi(k,s) = A(s)A(s+1)\cdots A(k-1)A(k)\qquad \hbox{for all $s$ and
$k$ with $k\ge s$,}
\]
where
\[
\Phi(k,k) = A(k)\qquad \hbox{for all $k$},
\]
and each $A(s)$ is a matrix whose $i$-th column is
the vector $a^i(s)$.
Using these transition matrices and the decomposition of
the estimate evolution of
Eqs.\ (\ref{decomp_dyn})--(\ref{nonlin_error}),
the relation between $x^i(k+1)$ and
the estimates $x^1(s),\ldots,x^m(s)$ at time $s\le k$ is given by
\begin{equation}
x^i(k+1) = \sum_{j=1}^m [\Phi(k,s)]^i_j \, x^j(s) +\sum_{r=s+1}^k
\left(\sum_{j=1}^m [\Phi(k,r)]^i_j\, e^j(r-1) \right) +e^i(k).
\label{short0}
\end{equation}
Here we can view $e^j(k)$ as an external perturbation input to the
system.

We use this relation to study the ``steady-state'' behavior of a
related process. In particular, we define an auxiliary sequence
$\{y(k)\}$, where $y(k)$ is given by
\begin{equation}
y(k)={1\over m} \sum_{i=1}^m w^i(k) \qquad\hbox{for all }k.
\label{defykwj}
\end{equation}
Since $w^i(k)=\sum_{j=1}^m a_j^i(k)x^j(k)$, under the doubly
stochasticity of the weights, it follows that
\begin{equation}
y(k)={1\over m} \sum_{j=1}^m x^j(k) \qquad\hbox{for all }k.
\label{defyk}
\end{equation}
Furthermore, from the relations in (\ref{short0}) using the doubly
stochasticity of the weights, we have
\begin{equation}
y(k)= {1\over m} \sum_{j=1}^m x^j(s) + {1\over m}\sum_{r=s+1}^{k}
\left(\sum_{j=1}^m  e^j(r-1)\right).\label{approxiterate}
\end{equation}

We now show that the limiting behavior of the agent estimates
$x^i(k)$ is the same as the limiting behavior of $y(k)$ as
$k\to\infty$. We establish this result using the assumptions on the
multi-agent model of Section \ref{sec:assumptions}.

\begin{lemma}
\emph{Let the intersection set $X=\cap_{i=1}^m X_i$ be nonempty.
Also, let Assumptions \ref{weightsrule}, \ref{doubstoc},
\ref{connected}, and \ref{boundedintervals} hold.
%\ref{weightsrule}, \ref{doubstoc}, \ref{connected}, and
%\ref{boundedintervals}].
We then have
\[\lim_{k\to \infty} \|x^i(k)-y(k)\|=0,\qquad \lim_{k\to \infty}
\|w^i(k)-y(k)\|=0\qquad \hbox{for all }i.\]
}\label{distwkyk}
\end{lemma}

\begin{proof}
By Lemma \ref{seqbehav}(d),
we have
$e^i(k)\to 0$ as $k\to \infty$ for all $i$.
Therefore, for any $\e>0$, we can choose some integer
$s$ such that $\|e^i(k)\|\le \e$ for all
$k\ge s$ and for all $i$.
Using the relations in Eqs.\ (\ref{short0}) and (\ref{approxiterate}),
we obtain for all $i$ and $k\ge s+1$,
\begin{eqnarray*}
\|x^i(k)& - & y(k)\|
= \left\|\sum_{j=1}^m \Big([\Phi(k-1,s)]_j^i -
{1\over m}\Big)x^j(s) \right.\\
& + & \left. \sum_{r=s+1}^{k-1} \sum_{j=1}^m \Big(
[\Phi(k-1,r)]_j^i - {1\over m}\Big)e^j(r-1)
+\Big(e^i(k-1)- {1\over m} \sum_{j=1}^m e^j(k-1)\Big)\right\|\\
&\le& \sum_{j=1}^m \Big|[\Phi(k-1,s)]_j^i - {1\over m}\Big|\
\|x^j(s)\| \\ & + & \sum_{r=s+1}^{k-1} \sum_{j=1}^m
\Big|[\Phi(k-1,r)]_j^i - {1\over m}\Big|
\|e^j(r-1)\| +  \|e^i(k-1)\| + {1\over m}\sum_{j=1}^m \|e^j(k-1)\|.\\
\end{eqnarray*}
Using the estimates for $\Big|[\Phi(k-1,s)]_j^i-{1\over m} \Big|$
of Proposition \ref{tranconv}(a), we have
\begin{eqnarray}
\|x^i(k)-y(k)\| &\le& 2\,{1+\eta^{-B_0}\over 1-\eta^{B_0}}\ \left(1-\eta^{B_0}
\right)^{{k-1-s\over B_0}}\ \sum_{j=1}^m \|x^j(s)\| \cr
&& +
\sum_{r=s+1}^{k-1}  2\,{1+\eta^{-B_0}\over 1-\eta^{B_0}}\
\left(1-\eta^{B_0}\right)^{{k-1-r\over B_0}} \,
\sum_{j=1}^m
\|e^j(r-1)\| \cr
&& + \,\|e^i(k-1)\| +{1\over m}\sum_{j=1}^m \|e^j(k-1)\|.
\label{xiandyk}
\end{eqnarray}
Since $\|e^i(k)\|\le \e$ for all $k\ge s$ and for all
$i$, from the preceding inequality we obtain
\begin{eqnarray*}
\|x^i(k)-y(k)\| &\le& 2\,{1+\eta^{-B_0}\over 1-\eta^{B_0}}\
\left(1-\eta^{B_0}\right)^{{k-1-s\over B_0}} \sum_{j=1}^m \|x^j(s)\| \\ && +
2 m \e\, {1+\eta^{-B_0}\over 1-\eta^{B_0}}\ {1\over 1-
(1-\eta^{B_0})^{{1\over B_0}}}  + 2\e.
\end{eqnarray*}
Thus, by taking the limit superior as $k\to\infty$, we see that
\[\limsup_{k\to\infty} \|x^i(k)-y(k)\|
\le 2 m \e\, {1+\eta^{-B_0}\over 1-\eta^{B_0}}\ {1\over 1-
(1-\eta^{B_0})^{{1\over B_0}}}  + 2\e,$$
which by the arbitrary choice of $\e$, implies
$\lim_{k\to \infty} \|x^i(k)-y(k)\|=0$
for all $i$.

Consider now $\sum_{i=1}^m\|w^i(k)-y(k)\|$.
By using $w^i(k)=\sum_{j=1}^m a^i_j(k)x^j(k)$ [cf.\ Eq.\ (\ref{wvector})]
and the stochasticity of the vector $a^i(k)$,
we have
$$\sum_{i=1}^m\|w^i(k)-y(k)\|\le
\sum_{i=1}^m\sum_{j=1}^m a^i_j(k)\|x^j(k) -y(k)\|.$$ By exchanging
the order of the summations over $i$ and $j$, and using the doubly
stochasticity of $a^i(k)$, we further obtain
\begin{equation}\label{sumwiandyk}
\sum_{i=1}^m\|w^i(k)-y(k)\|\le
\sum_{j=1}^m\left(\sum_{i=1}^m a^i_j(k)\right)\|x^j(k) -y(k)\|
=
\sum_{j=1}^m\|x^j(k) -y(k)\|.
\end{equation}
Since $\lim_{k\to \infty} \|x^j(k)-y(k)\|=0$
for all $j$, we have
$$\lim_{k\to \infty} \sum_{i=1}^m\|w^i(k)-y(k)\|=0,$$
implying $\lim_{k\to\infty} \|w^i(k)-y(k)\|=0$ for all $i$.
\end{proof}

\vskip 1pc

We next show that the agents reach a consensus asymptotically, i.e.,
the agent estimates $x^i(k)$ converge to the same point as $k$ goes
to infinity.

\begin{proposition}
\ {\it (Consensus)}\quad \emph{Let the set $X=\cap_{i=1}^m X_i$ be
nonempty. Also, let Assumptions \ref{weightsrule}, \ref{doubstoc},
\ref{connected}, and \ref{boundedintervals} hold.
%\ref{weightsrule}, \ref{doubstoc}, \ref{connected}, and
%\ref{boundedintervals}].
For all $i$, let the sequence $\{x^i(k)\}$
be generated by the projected consensus algorithm (\ref{proj-estimate}).
We then have for some $\tl x\in X$,
$$
\lim_{k\to\infty}\|x^i(k)-\tl x\|=0\qquad
\lim_{k\to\infty}\|w^i(k)-\tl x\|=0,\qquad
\hbox{for all }i.$$
}
\label{consensus}
\end{proposition}

\begin{proof}
The proof idea is to consider the sequence $\{y(k)\}$, defined in
Eq.\ (\ref{approxiterate}), and show that it has  a limit point in
the set $X$. By using this and Lemma
\ref{distwkyk}, we establish the convergence of
each $w^i(k)$ and $x^i(k)$ to $\tl x$.

To show that $\{y(k)\}$ has a limit point in the set $X$, we first
consider the sequence
$$\sum_{j=1}^m \dist(y(k),X_j).$$
Since $x^j(k)\in X_j$ for all $j$ and $k\ge 0$, we have
\[\sum_{j=1}^m \dist(y(k),X_j) \le \sum_{j=1}^m \|y(k)-x^j(k)\|.\]
Taking the limit as $k\to \infty$ in the preceding relation and
using Lemma \ref{distwkyk}, we conclude
\begin{equation}\lim_{k\to\infty}\
\sum_{j=1}^m\dist(y(k),X_j)=0.\label{zerodist}
\end{equation}
For a given $x\in X$, using Lemma \ref{seqbehav}(c), we have
\[\sum_{i=1}^m \|x^i(k) - x\| \le \sum_{i=1}^m \|x^i(0) - x\|\qquad \hbox{for all }k\ge 0.\]
This implies that the sequence $\{\sum_{i=1}^m \|x^i(k) - x\|\}$,
and therefore each of the sequences $\{x^i(k)\}$ are bounded. Since
for all $i$
\[\|y(k)\|\le \|x^i(k) - y(k)\|+\|x^i(k)\|\qquad \hbox{for all }k\ge 0,\]
using Lemma \ref{distwkyk}, it follows that the sequence $\{y(k)\}$
is bounded. In view of Eq.\ (\ref{zerodist}), this implies that the
sequence $\{y(k)\}$ has a limit point $\tl x$ that belongs to the
set $X=\cap_{j=1}^m X_j$. Furthermore, because
$\lim_{k\to\infty}\|w^i(k)-y(k)\|=0$ for all $i$, we conclude that
$\tl x$ is also a limit point of the sequence $\{w^i(k)\}$  for all
$i$. Since the sum sequence $\Big\{\sum_{i=1}^m \|w^i(k)-\tl
x\|\Big\}$ is nonincreasing by Lemma \ref{seqbehav}(c) and since
each $\{w^i(k)\}$ is converging to $\tl x$ along a subsequence, it
follows that
$$\lim_{k\to\infty}\ \sum_{i=1}^m \|w^i(k)-\tl x\|=0,$$
implying $\lim_{k\to \infty} \|w^i(k)-\tl x\|=0$ for all $i$.
Using this, together with the relations
$\lim_{k\to \infty} \|w^i(k)-y(k)\|=0$ and
$\lim_{k\to \infty} \|x^i(k)-y(k)\|=0$ for all $i$
(cf.\ Lemma \ref{distwkyk}),
we conclude
$$\lim_{k\to\infty} \|x^i(k)-\tl x\|=0
\qquad\hbox{for all }i.$$
\end{proof}

\subsection{Convergence Rate}\label{rate_analysis} In this section, we
establish a convergence rate result for the iterates $x^i(k)$
generated by the projected consensus algorithm (\ref{proj-estimate})
for the case when the weights are time-invariant and equal, i.e.,
$a^i(k)=(1/m,\ldots,1/m)'$ for all $i$ and $k$. In our multi-agent
model, this case corresponds to a fixed and complete connectivity
graph, where each agent is connected to every other agent. We
provide our rate estimate under an interior point assumption on the
sets $X_i$, stated in Assumption \ref{assump:intcond}.

We first establish a bound on the distance from the vectors of a
convergent sequence to the limit point of the sequence. This
relation holds for constant uniform weights, and it is motivated by
a similar estimate used in the analysis of alternating projections
methods in Gubin {\it et al.} \cite{gubinpolyak} (see the proof of
Lemma 6 there).

\begin{lemma} \label{wlim}\emph{
Let $Y$ be a nonempty closed convex set in $\rn$. Let
$\{u(k)\}\subseteq\rn$ be a sequence converging to some $\tilde y\in
Y$, and such that $\|u(k+1)-y\|\le \|u(k)-y\|$ for all $y\in Y$ and
all $k$. We then have \[ \|u(k)-\tilde y\|\le 2\,
\dist(u(k),Y)\qquad\hbox{for all }k\ge0.\] }
\end{lemma}

\begin{proof}
Let $B(x,\a)$ denote the closed ball centered at a vector $x$ with radius
$\a$, i.e., $B(x,\a) = \{z\ |\ \|z-x\|\le \a\}$. For each $l$,
consider the sets
\[S_l = \bigcap_{k=0}^l B\Big(P_Y[u(k)], \dist(u(k),Y)\Big).\]
The sets $S_l$ are convex, compact, and nested, i.e.,
$S_{l+1}\subseteq S_l$ for all $l$. The nonincreasing property
of the sequence $\{u(k)\}$ implies that $\|u(k+s)-P_Y[u(k)]\|\le
\|u(k)-P_Y[u(k)]\|$ for all $k,s\ge 0$; hence,
the sets $S_l$ are also nonempty. Consequently,
their intersection
$\cap_{l=0}^\infty S_l$ is nonempty and every point $y^*\in
\cap_{l=0}^\infty S_l$ is a limit point of the sequence $\{u(k)\}$.
By assumption, the sequence  $\{u(k)\}$ converges to $\tl y\in Y$,
and therefore, $\cap_{l=0}^\infty S_l=\{\tl y\}$. Then, in view of
the definition of the sets $S_l$, we obtain for all $k$,
\[\|u(k)-\tl y\|
\le \|u(k)-P_Y[u(k)]\| + \|P_Y[u(k)]-\tl y\|\le 2\,
\dist(u(k),Y).\]
\end{proof}

\vskip 1pc

We now establish a convergence rate result for constant uniform
weights. In  particular, we show that the projected consensus
algorithm converges with a geometric rate under the Interior Point
assumption.

\begin{proposition}
\emph{Let Assumptions \ref{assump:intcond}, \ref{weightsrule},
\ref{doubstoc}, \ref{connected}, and \ref{boundedintervals} hold.
Let the weight vectors $a^i(k)$ in algorithm (\ref{proj-estimate})
be given by $a^i(k)=(1/m ,\ldots,1/m)'$ for all $i$ and $k$. For all
$i$, let the sequence $\{x^i(k)\}$ be generated by the algorithm
(\ref{proj-estimate}). We then have
\[\sum_{i=1}^m \|x^i(k)-\tilde x\|^2
\le \left(1-\frac{1}{ 4 R^2}\right)^k \ \sum_{i=1}^m \|x^i(0)-\tilde
x\|^2 \qquad\hbox{for all }k\ge0,\] where $\tilde x\in X$ is the
limit of the sequence $\{x^i(k)\}$, and
$R=\frac{1}{\d}\,\sum_{i=1}^m\|x^i(0)-\bar x\|$ with $\bar x$ and
$\d$ given in the Interior Point assumption.}
\end{proposition}

\begin{proof}
Since the weight vectors $a^i(k)$ are given by
$a^i(k)=(1/m,\ldots,1/m)'$, it follows that
\[w^i(k)=w(k) = {1\over m} \sum_{j=1}^m x^j(k)\qquad \hbox{for all }i,\]
[see the definition of $w^i(k)$ in Eq.\ (\ref{wvector})]. For all
$k\ge 0$, using Lemma \ref{distance}(b) with the identification
$x^i=x^i(k)$ for each $i=1,\ldots,m$, and $\hat{x}=w(k)$, we obtain
\[ \dist(w(k),X) \le {1\over \delta m}\, \Big(\sum_{j=1}^m \|x^j (k) -\bar
x\|\Big) \Big(\sum_{j=1}^m \dist(w(k), X_j)\Big),\] where the vector
$\bar x$ and the scalar $\delta$ are given in Assumption
\ref{assump:intcond}. Since $\bar x\in X$, the sequence
$\{\sum_{i=1}^m \|x^i(k)-\bar x\|\}$ is nonincreasing by Lemma
\ref{seqbehav}(c). Therefore, we have $\sum_{i=1}^m \|x^i(k+1)-\bar
x\|\le \sum_{i=1}^m \|x^i(0)-\bar x\|$ for all $k$. Defining the
constant $R= {1\over \delta} \sum_{i=1}^m\|x^i(0)-\bar x\|$ and
substituting in the preceding relation, we obtain
\begin{eqnarray}
\dist(w(k),X) &\le& {R\over m}\,
\Big(\sum_{j=1}^m \dist(w(k), X_j)\Big)\nonumber\\
&=& {R\over m}\, \sum_{j=1}^m \|w(k)-x^j(k+1)\|,\label{raterel}
\end{eqnarray}
where the second relation follows in view of the definition of
$x^j(k+1)$ [cf.\ Eq.\ (\ref{error-estimate})].

By Proposition \ref{consensus}, we have $w(k)\to \tilde x$ for some
$\tilde x \in X$ as $k\to \infty$. Furthermore, by Lemma
\ref{seqbehav}(c) and the relation $w^i(k)=w(k)$ for all $i$ and
$k$, we have that the sequence $\{\|w(k)-x\|\}$ is nonincreasing for
any $x\in X$. Therefore, the sequence $\{w(k)\}$ satisfies the
conditions of Lemma \ref{wlim}, and by using this lemma we obtain
\[\|w(k)-\tilde x\|\le 2\, \dist(w(k),X)\qquad\hbox{for all }k.\]
Combining this relation with Eq.\ (\ref{raterel}), we further obtain
\[\|w(k)-\tilde x\|
\le {2R\over m}\, \sum_{i=1}^m \|w(k)-x^i(k+1)\|.\]
Taking the square
of both sides and using the convexity of the square function
$(\cdot)^2$, we have
\begin{equation}
\|w(k)-\tilde x\| ^2 \le {4R^2\over m}\, \sum_{i=1}^m
\|w(k)-x^i(k+1)\|^2.\label{yeter}\end{equation} Since
$x^i(k+1)=P_{X_i}[w(k)]$ for all $i$ and $k$, using Lemma
\ref{seqbehav}(a) with the substitutions $x=\tl x\in X$ and
$e^i(k)=x^i(k+1)-w(k)$ for all $i$, we see that
\[\sum_{i=1}^m \|w(k)-x^i(k+1)\|^2
\le m\,\|w(k)-\tl x\|^2- \sum_{i=1}^m\|x^i(k+1)-\tl x\|^2
\qquad \hbox{for all }k.\]
Using this relation in Eq.\
(\ref{yeter}), we obtain
\[
\|w(k)-\tilde x\| ^2 \le {4 R^2\over m}\,
\left( m\,\|w(k)-\tilde x\|^2 -
\sum_{i=1}^m\|x^i(k+1)-\tilde x\|^2\right).
\]
Rearranging the terms and using the relation
$m\,\|w(k)-\tilde x\|^2 \le \sum_{i=1}^m\|x^i(k)-\tilde
x\|^2$ [cf.\ Lemma \ref{seqbehav}(a) with $w(k)=w^i(k)$ and $x=\tl x$],
we obtain
\[\sum_{i=1}^m\|x^i(k+1)-\tilde x\|^2
\le \left(1-\frac{1}{4R^2}\right)\sum_{i=1}^m\|x^i(k)-\tilde
x\|^2,\] which yields the desired result.
\end{proof}

\section{Constrained Optimization}\label{sec:optim}

We next consider the problem of optimizing the sum of convex
objective functions corresponding to $m$ agents connected over a
time-varying topology. The goal of the agents is to cooperatively
solve the constrained optimization problem
\begin{eqnarray}
\hbox{minimize } && \sum_{i=1}^m f_i(x)\label{constoptim}\\
\hbox{subject to} && x\in \bigcap_{i=1}^m X_i,
\end{eqnarray}
where each $f_i:\re^n\to\re$ is a convex function, representing the
local objective function of agent $i$, and each $X_i\subseteq \rn$
is a closed convex set, representing the local constraint set of
agent $i$. We assume that the local objective function $f_i$ and the
local constraint set $X_i$ are known to agent $i$ only.

To keep our discussion general, we do not assume differentiability
of any of the functions $f_i$. Since each $f_i$ is convex over the
entire $\re^n$, the function is differentiable almost everywhere
(see \cite{book} or \cite{Rock70}). At the points where the function
fails to be differentiable, a subgradient exists and can be used in
the role of a gradient. In particular, for a given convex function
$F:\re^n\to\re$ and a point $\ol x$, {\it a subgradient of the
function $F$ at $\ol x$} is a vector $s_F(\ol x)\in\re^n$ such that
\begin{equation}
F(\ol x) + s_F(\ol x)'(x-\ol x)\le F(x) \qquad \hbox{for all }x.
\label{sgdconvdef}
\end{equation}
The set of all subgradients of $F$ at a given point $\bar x$ is
denoted by $\partial F(\bar x)$, and it is referred to as the {\it
subdifferential set} of $F$ at $\bar x$.

\subsection{Distributed Projected Subgradient Algorithm}

We introduce a distributed subgradient method for solving problem
(\ref{constoptim}) using the assumptions imposed on the information
exchange among the agents in Section~\ref{sec:assumptions}. The main
idea of the algorithm is the use of consensus as a mechanism for
distributing the computations among the agents. In particular, each
agent $i$ starts with an initial estimate $x^i(0)\in X_i$ and
updates its estimate. An agent $i$ updates its estimate by combining
the estimates received from its neighbors, by taking a subgradient
step to minimize its objective function $f_i$, and by projecting on
its constraint set $X_i$. Formally, each agent $i$ updates according
to the following rule:
\begin{eqnarray}
v^i(k)& = &\sum_{j=1}^m a^i_j(k)x^j(k)\label{combine}\\ x^i(k+1) & =
& P_{X_i}\left[v^i(k) - \a_k d_i(k)\right], \label{sgdstep}
\end{eqnarray}
where the scalars $a^i_j(k)$ are nonnegative weights and the scalar
$\a_k>0$ is a stepsize. The vector $d_i(k)$ is a subgradient of the
agent $i$ local objective function $f_i(x)$ at $x=v^i(k)$.

We refer to the method (\ref{combine})-(\ref{sgdstep}) as the {\it
projected subgradient algorithm}. To analyze this algorithm, we find
it convenient to
re-write the relation for $x^i(k+1)$  %(\ref{sgdstep})
in an equivalent form. This form helps us identify the linear
effects due to agents mixing the estimates [which will be driven by
the transition matrices $\Phi(k,s)$], and the nonlinear effects due
to taking subgradient steps and projecting. In particular, we
re-write the relations (\ref{combine})--(\ref{sgdstep}) as follows:
\begin{eqnarray}
v^i(k)  &= &\sum_{j=1}^m a^i_j(k)x^j(k)\cr x^i(k+1) &=& v^i(k)-\a_k
d_i(k)+\phi^i(k)\label{rewrite}\\ \phi^i(k) &=& P_{X_i}\left[v^i(k)
- \a_k d_i(k)\right] -\left(v^i(k)-\a_k d_i(k)\right).
\label{error_form}
\end{eqnarray}

The evolution of the iterates is complicated due to the nonlinear
effects of the projection operation, and even more complicated due
to the projections on different sets. In our subsequent analysis, we
study two special cases: 1) when the constraint sets are the same
[i.e., $X_i=X$ for all $i$], but the agent connectivity is
time-varying; and 2) when the constraint sets $X_i$ are different,
but the agent communication graph is fully connected.
In the analysis of both cases, we use a basic relation for the iterates
$x^i(k)$ generated by the method in (\ref{error_form}). The relation
is established in the following lemma.

\begin{lemma}\label{lemma:basic}\emph{
Let Assumptions \ref{weightsrule} and \ref{doubstoc} and hold. Let
$\{x^i(k)\}$ be the iterates generated by the algorithm
(\ref{combine})-(\ref{sgdstep}). We have for any $z\in
X=\cap_{i=1}^m X_i$ and all $k\ge0$,
\begin{eqnarray*}
\sum_{i=1}^m\|x^i(k+1)-z\|^2 &\le& \sum_{i=1}^m\|x^i(k)-z\|^2
+\a^2_k \sum_{i=1}^m \|d_i(k)\|^2
-2\a_k\sum_{i=1}^m\left(f_i(v^i(k))-f_i(z)\right)\\
&& -\sum_{i=1}^m\|\phi^i(k)\|^2.
\end{eqnarray*}}
\end{lemma}

\begin{proof}
Since  $x^i(k+1)=P_{X_i} [v^i(k)-\a_k d_i(k)],$ it follows from
Lemma \ref{proj_prop}(b) and from the definition of the projection
error $\phi^i(k)$ in (\ref{error_form}) that for any $z\in X,$
\[\|x^i(k+1)-z\|^2\le \|v^i(k)-\a_k d_i(k)-z\|^2-\|\phi^i(k)\|^2.\]
By expanding the term $\|v^i(k)-\a_k d_i(k)-z\|^2$, we obtain
\[\|v^i(k)-\a_k d_i(k)-z\|^2=\|v^i(k)-z\|^2 +\a^2_k \|d_i(k)\|^2
-2\a_k d_i(k)'(v^i(k)-z).\]
Since $d_i(k)$ is a subgradient of $f_i(x)$ at $x=v^i(k)$, we have
\[d_i(k)'(v^i(k)-z)\ge f_i(v^i(k))-f_i(z).\]
By combining the preceding relations, we obtain
\[\|x^i(k+1)-z\|^2\le \|v^i(k)-z\|^2+\a^2_k \|d_i(k)\|^2
-2\a_k\left(f_i(v^i(k))-f_i(z)\right)-\|\phi^i(k)\|^2.
\]

Since $v^i(k)=\sum_{j=1}^m a^i_j(k) x^j(k)$, using the convexity of
the norm square function and the stochasticity of the weights
$a_j^i(k)$, $j=1, \ldots,m$, it follows that
\[\|v^i(k)-z\|^2\le \sum_{j=1}^m a^i_j(k) \|x^j(k)-z\|^2.\]
Combining the preceding two relations, we obtain
\begin{eqnarray*}
\|x^i(k+1)-z\|^2 &\le&  \sum_{j=1}^m a^i_j(k) \|x^j(k)-z\|^2+\a^2_k
\|d_i(k)\|^2
-2\a_k\left(f_i(v^i(k))-f_i(z)\right)\\
&& -\|\phi^i(k)\|^2.
\end{eqnarray*}
By summing the preceding relation over $i=1,\ldots,m,$ and using the
doubly stochasticity of the weights, i.e.,
\[\sum_{i=1}^m \sum_{j=1}^m a^i_j(k) \|x^j(k)-z\|^2
=\sum_{j=1}^m\left(\sum_{i=1}^m a^i_j(k)\right) \|x^j(k)-z\|^2
=\sum_{j=1}^m \|x^j(k)-z\|^2,\]
we obtain the desired relation.
\end{proof}

\subsubsection{Convergence when $X_i=X$ for all $i$}

We first study the case when all constraint sets are the same, i.e.,
$X_i=X$ for all $i$. The next assumption formally states the
conditions we adopt in the convergence analysis.

\begin{assumption} \emph{
\begin{itemize}
\item [(a)]
The constraint sets $X_i$ are the same, i.e, $X_i=X$ for a closed
convex set $X$.
\item[(b)]
The subgradient sets of each $f_i$ are bounded over the set $X$,
i.e., there is a scalar $L>0$ such that for all $i$,
\[\|d\|\le L\qquad\hbox{for all }d\in\partial f_i(x)
\quad\hbox{and all }x\in X.\]
\end{itemize}\label{assump:equal_sets}}
\end{assumption}
The subgradient boundedness assumption in part (b) holds for example
when the set $X$ is compact (see \cite{book}).

In proving our convergence results, we use a property of the
infinite sum of products of the components of two sequences. In
particular, for a scalar $\beta\in (0,1)$ and a scalar sequence
$\{\gamma_k\}$, we consider the ``convolution'' sequence
$\sum_{\ell=0}^k\beta^{k-\ell}\gamma_\ell =\beta^k\gamma_0
+\beta^{k-1}\gamma_1+\cdots+\beta\gamma_{k-1}+\gamma_k.$ We have the
following result.
\begin{lemma} \label{lemma:seq}\emph{
Let $0< \beta<1$ and let  $\{\gamma_k\}$ be a positive scalar
sequence. Assume that $\lim_{k\to\infty}\gamma_k=0.$ Then
\[\lim_{k\to\infty}\sum_{\ell=0}^k \beta^{k-\ell}\gamma_\ell=0.\]
In addition, if $\sum_{k}\gamma_k<\infty,$ then
\[\sum_{k}\sum_{\ell=0}^k \beta^{k-\ell}\gamma_\ell<\infty.\]}
\end{lemma}

\begin{proof}
Let $\epsilon>0$ be arbitrary.  Since $\a_k\to0$, there is an index
$K$ such that $\a_k\le\epsilon$ for all $k\ge K$. For all $k\ge
K+1$, we have
\[\sum_{\ell=0}^k\beta^{k-\ell}\gamma_\ell=
\sum_{\ell=0}^K\beta^{k-\ell}\gamma_\ell +
\sum_{\ell=K+1}^k\beta^{k-\ell}\gamma_\ell \le \max_{0\le t\le
K}\gamma_t \sum_{\ell=0}^K\beta^{k-\ell}
+\epsilon\sum_{\ell=K+1}^k\beta^{k-\ell}.\] Since
$\sum_{\ell=K+1}^k\beta^{k-\ell}\le\frac{1}{1-\beta}$ and
\[\sum_{\ell=0}^K\beta^{k-\ell}=\beta^k+\cdots+\beta^{k-K}
=\beta^{k-K}(1+\cdots+\beta^K)\le \frac{\beta^{k-K}}{1-\beta},\] it
follows that for all $k\ge K+1$,
\[\sum_{\ell=0}^k\beta^{k-\ell}\gamma_\ell
\le \max_{0\le t\le K}\gamma_t\,\frac{\beta^{k-K}}{1-\beta}
+\frac{\epsilon}{1-\beta}.\] Therefore,
\[\limsup_{k\to\infty}\sum_{\ell=0}^k\beta^{k-\ell}\gamma_\ell
\le \frac{\epsilon}{1-\beta}.\] Since $\epsilon$ is arbitrary, we
conclude that
$\limsup_{k\to\infty}\sum_{\ell=0}^k\beta^{k-\ell}\gamma_\ell=0$,
implying
\[\lim_{k\to\infty}\sum_{\ell=0}^k\beta^{k-\ell}\gamma_\ell=0.\]

Suppose now $\sum_k\gamma_k<\infty.$ Then, for any integer $M\ge1$,
we have
\[\sum_{k=0}^M \left(\sum_{\ell=0}^k\beta^{k-\ell}\gamma_\ell\right)
=\sum_{\ell=0}^M\gamma_\ell\sum_{t=0}^{M-\ell}\beta^t \le
\sum_{\ell=0}^M\gamma_\ell\frac{1}{1-\beta},\] implying that
\[\sum_{k=0}^\infty \left(\sum_{\ell=0}^k\beta^{k-\ell}\gamma_\ell\right)
\le \frac{1}{1-\beta}\sum_{\ell=0}^\infty\gamma_\ell<\infty.\]
\end{proof}

Our goal is to show that the agent disagreements $\|x^i(k)-x^j(k)\|$
converge to zero. To measure the agent disagreements
$\|x^i(k)-x^j(k)\|$ in time, we consider their average
$\frac{1}{m}\sum_{j=1}^m x^j(k)$, and consider the agent
disagreement with respect to this average. In particular, we define
\[y(k)=\frac{1}{m}\sum_{j=1}^m x^j(k)\qquad\hbox{for all }k.\]
In view of Eq.\ (\ref{rewrite}), we have
\[y(k+1)=\frac{1}{m}\sum_{i=1}^m v^i(k) -\frac{\a_k }{m}\sum_{i=1}^m d_i(k)
+\frac{1}{m}\sum_{i=1}^m \phi^i(k).\] When the weights are doubly
stochastic, since $v^i(k)=\sum_{j=1}^m a^i_j(k) x^j(k)$, it follows
that
\begin{equation}
y(k+1)=y(k) -\frac{\a_k }{m}\sum_{i=1}^m d_i(k)
+\frac{1}{m}\sum_{i=1}^m \phi^i(k). \label{y_evol}
\end{equation}
Under Assumption \ref{assump:equal_sets}, the assumptions on the
agent weights and connectivity stated in Section
\ref{sec:assumptions}, and some conditions on the stepsize $\a_k$,
the next lemma studies the convergence properties of the sequences
$\Big\{\|x^i(k)-y(k)\|\Big\}$ for all $i$.

\begin{lemma}\label{xkyktozero}\emph{
Let Assumptions~\ref{weightsrule}, \ref{doubstoc}, \ref{connected},
\ref{boundedintervals}, and \ref{assump:equal_sets} hold. Let
$\{x^i(k)\}$ be the iterates generated by the algorithm
(\ref{combine})-(\ref{sgdstep}) and consider the auxiliary sequence
$\{y(k)\}$ defined in (\ref{y_evol}).
\begin{itemize}
\item[(a)]
If the stepsize satisfies $\lim_{k\to\infty}\a_k=0$, then
\[\lim_{k\to\infty}\|x^i(k)-y(k)\|=0\qquad\hbox{for all }i.\]
\item[(b)]
If the stepsize satisfies $\sum_{k\to\infty}\a_k^2<\infty$, then
\[\sum_{k=1}^\infty\a_k
\|x^i(k)-y(k)\|<\infty\qquad\hbox{for all }i.\]
\end{itemize}}
\end{lemma}

\begin{proof}
(a) \ Using the relations in (\ref{error_form}) and the transition
matrices $\Phi(k,s)$, we can write for all $i$, and for all $k$ and
$s$ with $k>s,$
\begin{eqnarray*}
x^i(k+1) &=&\sum_{j=1}^m [\Phi(k,s)]^i_j x^j(s)- \sum_{r=s}^{k-1}
\sum_{j=1}^m [\Phi(k,r+1)]^i_j \a_r  d_j(r) - \a_k  d_i(k)\cr
&& + \sum_{r=s}^{k-1} \sum_{j=1}^m [\Phi(k,r+1)]^i_j \phi^j(r) +
\phi^i(k).
%\label{trans}
\end{eqnarray*}
Similarly, using the transition matrices and relation
(\ref{y_evol}), we can write for $y(k+1)$ and for all $k$ and $s$
with $k>s,$
\begin{eqnarray*}
y(k+1)= y(s)- \frac{1}{m}\sum_{r=s}^{k-1} \sum_{j=1}^m \a_r
d_j(r) - \frac{\a_k}{m}\sum_{i=1}^m  d_i(k)
 + \frac{1}{m}\sum_{r=s}^{k-1} \sum_{j=1}^m \phi^j(r)
+ \frac{1}{m}\sum_{j=1}^m \phi^j(k).
%\label{yk_trans}
\end{eqnarray*}
Therefore, since $y(s)=\frac{1}{m} \sum_{j=1}^m x^j(s)$, we have for
$s=0,$
\begin{eqnarray*}
\|x^i(k)-y(k)\| &\le & \sum_{j=1}^m \left|[\Phi(k-1,0)]^i_j
-\frac{1}{m}\right|\,\|x^j(0)\| \cr &&+ \sum_{r=0}^{k-2}
\sum_{j=1}^m \left|[\Phi(k-1,r+1)]^i_j- \frac{1}{m}\right| \, \a_r
\| d_j(r)\|\cr &&+ \a_{k-1}\|d_i(k-1)\|
+\frac{\a_{k-1}}{m}\sum_{j=1}^m \| d_j(k-1)\|\cr && +
\sum_{r=0}^{k-2} \sum_{j=1}^m
\left|[\Phi(k-1,r+1)]^i_j-\frac{1}{m}\right| \| \phi^j(r)\| \cr &&+
\|\phi^i(k-1)\|+\frac{1}{m}\sum_{j=1}^m \|\phi^j(k-1)\|.
\end{eqnarray*}
Using the estimate for $\left|[\Phi(k,s)]^i_j -\frac{1}{m}\right|$
of Proposition~\ref{tranconv}, we have for all $k\ge s$,
\[\left|[\Phi(k,s)]^i_j -\frac{1}{m}\right|\le C\beta^{k-s}\qquad
\hbox{for all }i,j,\] with $C=2\frac{1+\eta^{-B_0}}{1-\eta^{B_0}}$
and $\beta=\left(1-\eta^{B_0}\right)^{\frac{1}{B_0}}.$ Hence, using
this relation and the subgradient boundedness, we obtain for all $i$
and $k\ge2,$
\begin{eqnarray}
\|x^i(k)-y(k)\| &\le & mC\beta^{k-1}\sum_{j=1}^m \|x^j(0)\| + mCL
\sum_{r=0}^{k-2} \beta^ {k-r}\a_r + 2\a_{k-1}L\cr && +
C\sum_{r=0}^{k-2} \beta^ {k-r}\sum_{j=1}^m \| \phi^j(r)\| +
\|\phi^i(k-1)\|+\frac{1}{m}\sum_{j=1}^m \|\phi^j(k-1)\|.
\qquad\label{est_xkyk}
\end{eqnarray}

We next show that the errors $\phi^i(k)$ satisfy $\|\phi^i(k)\|\le
\a_k L$ for all $i$ and $k.$ In view of the relations in
(\ref{error_form}), since $x^j(k)\in X_j=X$ for all $k$ and $j$, and
the vector $a^i(k)$ is stochastic for all $i$ and $k$, it follows
that $v^i(k)\in X$ for all $i$ and $k$. Furthermore, by the
projection property in Lemma~\ref{proj_prop}(b), we have for all $i$
and $k$,
\begin{eqnarray*}
\|x^i(k+1)-v^i(k)\|^2 &\le& \|v^i(k)-\a_k d_i(k)-v^i(k)\|^2-
\|x^i(k+1)-(v^i(k)-\a_k d_i(k))\|^2\cr &\le& \a^2_k L^2
-\|\phi^i(k)\|^2,
\end{eqnarray*}
where in the last inequality we use $\|d_i(k)\|\le L$ (see
Assumption~\ref{assump:equal_sets}). It follows that
$\|\phi^i(k)\|\le \a_k L$ for all $i$ and $k.$ By using this in
relation (\ref{est_xkyk}), we obtain
\begin{eqnarray}
\|x^i(k)-y(k)\| &\le & mC\beta^{k-1}\sum_{j=1}^m \|x^j(0)\| + 2mCL
\sum_{r=0}^{k-2} \beta^ {k-r}\a_r + 4\a_{k-1}L.
\qquad\label{est_xkykak}
\end{eqnarray}

By taking the limit superior in relation (\ref{est_xkykak}) and
using the facts $\b^k\to0$ (recall $0<\beta<1$) and $\a_k\to0$, we
obtain for all $i,$
\[
\limsup_{k\to\infty}\|x^i(k)-y(k)\| \le 2mCL  \limsup_{k\to\infty}\,
\sum_{r=0}^{k-2} \beta^ {k-r}\a_r
\]
Finally, since $0<\beta<1$ and $\lim_{k\to\infty}\a_k=0$, by
Lemma~\ref{lemma:seq} we have
\[\lim_{k\to\infty}
\sum_{r=0}^{k-2} \beta^ {k-r}\a_r=0.\] In view of the preceding two
relations, it follows that $\lim_{k\to\infty}  \|x^i(k)-y(k)\|=0$
for all $i$.

\vskip 1pc

\noindent(b) \ By multiplying the relation in (\ref{est_xkykak})
with $\a_k$, we obtain
\[
\a_k\|x^i(k)-y(k)\| \le mC\a_k \beta^{k-1}\sum_{j=1}^m \|x^j(0)\| +
2mCL \sum_{r=0}^{k-2} \beta^ {k-r}\a_k \a_r + 4\a_k \a_{k-1}L.
\]
By using $\a_k \beta^{k-1}\le \a_k^2 +\beta^{2(k-1)}$ and
$2\a_k\a_r\le \a_k^2+\a_r^2$ for any $k$ and $r$, we have
\[\a_k\|x^i(k)-y(k)\| \le
mC\beta^{2(k-1)}\sum_{j=1}^m \|x^j(0)\| + mCA\a_k^2 + mCL
\sum_{r=0}^{k-2} \beta^ {k-r}\a_r^2 + 2L(\a_k^2 +\a_{k-1}^2),
\]
where $A= \sum_{j=1}^m \|x^j(0)\| + {L \over (1-\beta)}$. Therefore,
by summing and grouping some of the terms, we obtain
\begin{eqnarray*}
\sum_{k=1}^\infty \a_k\|x^i(k)-y(k)\| &\le &
mC\left(\sum_{k=1}^\infty \beta^{2(k-1)}\right) \sum_{j=1}^m
\|x^j(0)\| \cr &+& \sum_{k=1}^\infty \left(mCA\a_k^2 + 2L(\a_k^2
+\a_{k-1}^2)\right) +mCL \sum_{k=1}^\infty  \sum_{r=0}^{k-2} \beta^
{k-r}\a_r^2.
\end{eqnarray*}
In the preceding relation, the first term is summable since
$0<\b<1$. The second term is summable since $\sum_{k}\a_k^2<\infty$.
The third term is also summable by Lemma~\ref{lemma:seq}. Hence,
$\sum_{k=1}^\infty \a_k\|x^i(k)-y(k)\| <\infty.$
\end{proof}

\vskip 1pc

Using Lemmas \ref{lemma:basic} and \ref{xkyktozero}, we next show
that the iterates $x^i(k)$ converge to an optimal solution when we
use a stepsize converging to zero fast enough.

\begin{proposition}\label{porp:converg}\emph{
Let Assumptions~\ref{weightsrule}, \ref{doubstoc}, \ref{connected},
\ref{boundedintervals}, and \ref{assump:equal_sets} hold. Let
$\{x^i(k)\}$ be the iterates generated by the algorithm
(\ref{combine})-(\ref{sgdstep}) with the stepsize satisfying
$\sum_{k}\a_k =\infty$ and $\sum_{k}\a_k^2<\infty.$ In addition,
assume that the optimal solution set $X^*$ is nonempty. Then, there
exists an optimal point $x^*\in X^*$ such that
\[\lim_{k\to\infty}\|x^i(k)-x^*\|=0\qquad\hbox{for all }i.\]}
\end{proposition}

\begin{proof}
From Lemma \ref{lemma:basic}, we have for $z\in X$ and all $k,$
\begin{eqnarray*}
\sum_{i=1}^m\|x^i(k+1)-z\|^2 &\le& \sum_{j=1}^m\|x^j(k)-z\|^2
+\a^2_k\sum_{i=1}^m\|d_i(k)\|^2 \cr
&&-2\a_k\sum_{i=1}^m\left(f_i(v^i(k))-f_i(z)\right)
-\sum_{i=1}^m\|\phi^i(k)\|^2.
\end{eqnarray*}
By dropping the nonpositive term on the right hand side,
and by using the subgradient boundedness, we obtain
\begin{eqnarray}
\sum_{i=1}^m\|x^i(k+1)-z\|^2
&\le& \sum_{j=1}^m\|x^j(k)-z\|^2
+\a^2_k mL^2
-2\a_k\sum_{i=1}^m\left(f_i(v^i(k))-f_i(y(k))\right)\cr
&&-2\a_k\left(f(y(k))-f(z)\right).\qquad \label{eqn:mid}
\end{eqnarray}
In view of the subgradient boundedness and the stochasticity of the
weights, it follows
\[|f_i(v^i(k))-f_i(y(k))|\le L\|v^i(k)-y(k)\|
%=\left\|\sum_{j=1}^m a^i_j(k) x^j(k)-y(k)\right\|
\le L \sum_{j=1}^m a^i_j(k) \|x^j(k)-y(k)\|,\] implying, by the
doubly stochasticity of the weights, that
\[\sum_{i=1}^m\left|f_i(v^i(k))-f_i(y(k))\right|\le
L \sum_{j=1}^m\left(\sum_{i=1}^m a^i_j(k)\right) \|x^j(k)-y(k)\| = L
\sum_{j=1}^m \|x^j(k)-y(k)\|.\] By using this in relation
(\ref{eqn:mid}), we see that for any $z\in X$, and all $i$ and $k,$
\begin{eqnarray*}
\sum_{i=1}^m\|x^i(k+1)-z\|^2 &\le& \sum_{j=1}^m\|x^j(k)-z\|^2
+\a^2_k mL^2 +2\a_k L\sum_{j=1}^m \|x^j(k)- y(k)\|\cr
&&-2\a_k\left(f(y(k))-f(z)\right).
\end{eqnarray*}
By letting $z=z^*\in X^*$, and by re-arranging the terms and summing
these relations over some arbitrary window from $K$ to $N$ with
$K<N$, we obtain for any $z^*\in X^*$,
\begin{eqnarray}
\sum_{i=1}^m\|x^i(N+1)-z^*\|^2 &+&2\sum_{k=K}^{N}\a_k\, \left(
f(y(k))-f(z^*)\right) \le \sum_{i=1}^m\|x^i(K)-z^*\|^2\cr
&&+mL^2\sum_{k=K}^{N}\a^2_k
+2L \sum_{k=K}^{N} \a_k\sum_{j=1}^m
\|x^j(k)- y(k)\|. \label{eqn:one}
\end{eqnarray}

By letting $K=1$ and  $N\to\infty$ in relation (\ref{eqn:one}), and
using $\sum_{k=1}^{\infty}\a^2_k<\infty$ and $\sum_{k=1}^{\infty}
\a_k\sum_{j=1}^m \|x^j(k)- y(k)\|<\infty$ [which follows by
Lemma~\ref{xkyktozero}], we obtain
\[
\sum_{k=1}^{\infty}\a_k\, \left( f(y(k))-f(z^*)\right)<\infty.
\]
Since $x^j(k)\in X$ for all $j$, we have $y(k)\in X$ for all $k$.
Since $z^*\in X^*$, it follows that $f(y(k))-f^*\ge0$ for all $k$.
This relation, the assumption that $\sum_{k=1}^{\infty}\a_k=\infty$,
and $\sum_{k=1}^{\infty}\a_k\, \left( f(y(k))-f(z^*)\right)<\infty$
imply
\begin{equation}
\liminf_{k\to\infty} \left( f(y(k))-f^*\right)=0.
\label{eqn:liminffk}
\end{equation}

We next show that each sequence $\{x^i(k)\}$ converges to the same
optimal point. By dropping the nonnegative term involving
$f(y(k))-f(z^*)$ in (\ref{eqn:one}), we have
\begin{eqnarray*}
\sum_{i=1}^m\|x^i(N+1)-z^*\|^2 \le \sum_{i=1}^m\|x^i(K)-z^*\|^2
+mL^2\sum_{k=K}^{N}\a^2_k +2L \sum_{k=K}^{N} \a_k\sum_{j=1}^m
\|x^j(k)- y(k)\|.
\end{eqnarray*}
Since $\sum_{k}\a_k^2<\infty$ and $\sum_{k=1}^{\infty}
\a_k\sum_{j=1}^m \|x^j(k)- y(k)\|<\infty$, it follows that the
sequence $\{x^i(k)\}$ is bounded for each $i,$ and
\[\limsup_{N\to\infty}\,
\sum_{i=1}^m\|x^i(N+1)-z^*\|^2 \le \liminf_{K\to\infty}
\sum_{i=1}^m\|x^i(K)-z^*\|^2 \qquad\hbox{for all }i.\] Thus, the
scalar sequence $\{\sum_{i=1}^m\|x^i(k)-z^*\|\}$ is convergent for
every $z^*\in X^*$. By Lemma~\ref{xkyktozero}, we have
$\lim_{k\to\infty}\|x^i(k)-y(k)\|=0$. Therefore, it also follows
that $\{y(k)\}$ is bounded and the scalar sequence
$\{\|y(k)-z^*\|\}$ is convergent for every $z^*\in X^*$. Since
$y(k)$ is bounded, it must have a limit point, and in view of
$\liminf_{k\to\infty} f(y(k))=f^*$ [cf.\ Eq.\ (\ref{eqn:liminffk})]
and the continuity of $f$ (due to convexity of $f$ over $\rn$), one
of the limit points of $\{y(k)\}$ must belong to $X^*$; denote this
limit point by $x^*$. Since the sequence $\{\|y(k)-x^*\|\}$ is
convergent, it follows that $y(k)$ can have a unique limit point,
i.e., $\lim_{k\to \infty} y(k)=x^*$. This and
$\lim_{k\to\infty}\|x^i(k)-y(k)\|=0$ imply that each of the
sequences $\{x^i(k)\}$ converges to the same $x^*\in X^*$.
\end{proof}

\subsubsection{Convergence for uniform weights}
We next consider a version of the projected subgradient algorithm
(\ref{combine})--(\ref{sgdstep}) for the case when the agents use
uniform weights, i.e., $a_j^i(k)=\frac{1}{m}$ for all $i$, $j$, and
$k\ge 0$. We show that the estimates generated by the method
converge to an optimal solution of problem (\ref{constoptim}) under
some conditions. In particular, we adopt the following assumption in
our analysis.

\begin{assumption} \emph{For each $i$, the local constraint set $X_i$ is a compact set, i.e., there exists a scalar
$B>0$ such that
\[\|x\|\le B\qquad \hbox{for all }x\in X_i \hbox{ and all }i.\]
} \label{assump:compact}
\end{assumption}

An important implication of the preceding assumption is
that, for each $i$, the subgradients of the function $f_i$ at all
points $x\in X_i$ are uniformly bounded, i.e., there exists a scalar
$L>0$ such that
\begin{equation}
\|g\|\le L\qquad \hbox{for all }g\in \partial f_i(x), \hbox{ all
}x\in X_i \hbox{ and all }i.\label{sgdbd}
\end{equation}

Under this and the interior point assumption on the intersection set
$X=\cap_{i=1}^m X_i$ (cf.\ Assumption \ref{assump:intcond}), we have
the following result.

\begin{proposition} \emph{Let Assumptions
\ref{assump:intcond} and \ref{assump:compact} and hold. Let
$\{x^i(k)\}$ be the iterates generated by the algorithm
(\ref{combine})-(\ref{sgdstep}) with the weight vectors
$a^i(k)=(1/m ,\ldots,1/m)'$ for all $i$ and $k$, and the
stepsize satisfying $\sum_k \a_k = \infty$ and $\sum_k
\a_k^2<\infty$. Then, the sequences $\{x^i(k)\}$, $i=1,\ldots,m,$
converge to the same optimal point, i.e.,
\[\lim_{k\to\infty} x^i(k)=x^*\qquad
\hbox{for some $x^*\in X^*$ and all $i$}.\] \label{lemma:projzero}}
\end{proposition}

\begin{proof} By Assumption \ref{assump:compact}, each set $X_i$ is
compact, which implies that the intersection set $X=\cap_{i=1}^m
X_i$ is compact. Since each function $f_i$ is continuous (due to
being convex over $\rn$), it follows from Weierstrass' Theorem that
problem (\ref{constoptim}) has an optimal solution, denoted by
$z^*\in X$. By using Lemma \ref{lemma:basic} with $z=z^*$, we have
for all $i$ and $k\ge0,$
\begin{eqnarray}
\sum_{i=1}^m\|x^i(k+1)-z^*\|^2 &\le& \sum_{i=1}^m\|x^i(k)-z^*\|^2
+\a^2_k \sum_{i=1}^m \|d_i(k)\|^2\cr
&&-2\a_k\sum_{i=1}^m\left(f_i(v^i(k))-f_i(z^*)\right)
-\sum_{i=1}^m\|\phi^i(k)\|^2.\label{main-rel}
\end{eqnarray}
For any $k\ge 0$, define the vector $s(k)$ by
\[s(k)=\frac{\e}{\e+\delta}\, \bar x
+ \frac{\delta}{\e+\delta}\, \hat x(k),\] where $\hat x(k)={1\over
m} \sum_{i=1}^m x^i(k)$, $\e = \sum_{j=1}^m \dist(\hat x(k), X_j)$,
and $\delta$ is the scalar given in Assumption \ref{assump:intcond}
(cf.\ Lemma \ref{distance}). By using the subgradient boundedness
[see (\ref{sgdbd})] and adding and subtracting the term $2\a_k
\sum_{i=1}^m f_i(s(k))$ in Eq.\ (\ref{main-rel}), we obtain
\begin{eqnarray*}
\sum_{i=1}^m\|x^i(k+1)-z^*\|^2
&\le& \sum_{i=1}^m\|x^i(k)-z^*\|^2
+\a^2_k m L^2
-\sum_{i=1}^m\|\phi^i(k)\|^2\\ &&
-2\a_k\sum_{i=1}^m\left(f_i(s(k))-f_i(z^*)\right)
-2\a_k\sum_{i=1}^m\left(f_i(v^i(k))-f_i(s(k))\right).
\end{eqnarray*}
Using the subgradient definition and the subgradient boundedness
assumption, we further have
\[|f_i(v^i(k))-f_i(s(k))| \le L \|v^i(k)-s(k)\|\qquad
\hbox{for all $i$ and $k$}.\] Combining these relations with the
preceding and using the notation $f=\sum_{i=1}^m f_i$, we obtain
\begin{eqnarray}
\sum_{i=1}^m\|x^i(k+1)-z^*\|^2 &\le& \sum_{i=1}^m\|x^i(k)-z^*\|^2
+\a^2_k m L^2
-\sum_{i=1}^m\|\phi^i(k)\|^2\cr &&
-2\a_k\left(f(s(k))-f(z^*)\right)+2\a_k L \sum_{i=1}^m
\|v^i(k)-s(k)\|.\label{keymodrel}
\end{eqnarray}
Since the weights are all equal, from relation (\ref{combine})
we have $v^i(k)=\hat x(k)$ for all $i$ and $k$.
Using Lemma \ref{distance}(b) with the substitution $s=s(k)$ and
$\hat x=\hat x(k) = {1\over m}\sum_{j=1}^m x^j(k)$, we obtain
\[\|v^i(k)-s(k)\|\le \frac{1}{\delta m}
\Big(\sum_{j=1}^m \|x^j(k)-\bar x\|\Big)
\Big(\sum_{j=1}^m \dist(\hat x(k), X_j)\Big)\qquad\hbox{for all $i$ and
$k$}.\]
Since $x^j(k)\in X_j$, we have $\dist(\hat x(k), X_j)
\le \|\hat x(k)-x^j(k+1)\|$
for all $j$ and $k$, Furthermore, since
$\bar x \in X\subseteq X_j$ for all $j$,
using
Assumption~\ref{assump:compact}, we obtain $\|x^j(k)-\bar x\|\le 2B$.
Therefore, for all $i$ and
$k$,
\begin{equation}\|v^i(k)-s(k)\|
\le \frac{2B}{\delta}
\sum_{j=1}^m \dist(\hat x(k), X_j)
\le \frac{2B}{\delta}
\sum_{j=1}^m \|\hat x(k)-x^j(k+1)\|.
\label{pre-error}
\end{equation}
Moreover, we have $\hat x(k)=v^j(k)$ for all $j$ and $k$, implying
\[\|x^j(k+1)-\hat x(k)\|=\|x^j(k+1)-(v^j(k)-\a_kd_j(k))\|+\a_k \|d_j(k)\|.\]
In view of the definition of the error term $\phi^i(k)$ in (\ref{error_form})
and the subgradient boundedness, it follows
\[\|x^j(k+1)-\hat x(k)\|\le \|\phi^j(k)\|+\a_k L,\]
which when substituted in relation (\ref{pre-error}) yields
\begin{equation}
\|v^i(k)-s(k)\| \le \frac{2B}{\delta} \left(\a_k mL+ \sum_{j=1}^m
\|\phi^j(k)\|\right)\qquad\hbox{for all $i$ and $k$}.\label{projdis}
\end{equation}
We now substitute the estimate (\ref{projdis}) in Eq.\
(\ref{keymodrel}) and obtain for all $k$,
\begin{eqnarray}
\sum_{i=1}^m
\|x^i(k+1)-z^*\|^2
&\le& \sum_{i=1}^m\|x^i(k)-z^*\|^2
+\a^2_k m L^2 -\sum_{i=1}^m\|\phi^i(k)\|^2\cr
&& -2\a_k\left(f(s(k))-f(z^*)\right) +\frac{4m^2BL^2}{\delta}\,\a_k^2\cr
&&+{4\a_k m BL \over \delta} \sum_{i=1}^m
\|\phi^i(k)\|. \label{eqn:good}
\end{eqnarray}
Note that for each $i$, we can write
\begin{eqnarray*}
{4\a_k m BL \over \delta} \|\phi^i(k)\|
&=& 2
\left(\frac{2\sqrt{2}\a_k m BL}{\delta}\right)
\left(\frac{1}{\sqrt{2}}\|\phi^i(k)\|\right)\cr
&\le&\left(\frac{2\sqrt{2}\a_k m BL}{\delta}\right)^2 +
\frac{1}{2}\|\phi^i(k)\|^2.
\end{eqnarray*}
Therefore, by summing the preceding relations over $i$, we have for
all $k$,
\[{4\a_k m BL \over \delta}\sum_{i=1}^m \|\phi^i(k)\|
\le \frac{8m^3 B^2L^2}{\delta^2}\, \a_k^2 +\frac{1}{2}\sum_{i=1}^m
\|\phi^i(k)\|^2,\]
which when substituted in Eq.~(\ref{eqn:good}) yields
\begin{eqnarray*}
\sum_{i=1}^m\|x^i(k+1)-z^*\|^2 &\le& \sum_{i=1}^m\|x^i(k)-z^*\|^2
+ C\a_k^2
- \frac{1}{2}\sum_{i=1}^m\|\phi^i(k)\|^2\cr
&&-2\a_k\left(f(s(k))-f(z^*)\right),
\end{eqnarray*}
where $C=m L^2 +\frac{4m^2BL^2}{\delta}+\frac{8m^3 B^2L^2}{\delta^2}$.
By re-arranging the terms and summing the preceding relations over
$k$ for $k=K,\ldots, N$ for some arbitrary $K$ and $N$ with $K<N$,
we obtain
\begin{eqnarray}
\sum_{i=1}^m\|x^i(N+1)-z^*\|^2 &+&\frac{1}{2}\sum_{k=K}^N \sum_{i=1}^m
\|\phi^i(k)\|^2 +
2\sum_{k=K}^N\a_k\left(f(s(k))-f(z^*)\right)\cr &\le&
\sum_{i=1}^m\|x^i(K)-z^*\|^2 + C\sum_{k=K}^N\a_k^2. \label{eqn:ok}
\end{eqnarray}
By setting $K=0$ and letting $N\to\infty$, in view of
$\sum_{k}\a_k^2<\infty$, we see that
\[\frac{1}{2}\sum_{k=0}^\infty \sum_{i=1}^m
\|\phi^i(k)\|^2 +
2\sum_{k=0}^\infty\a_k\left(f(s(k))-f(z^*)\right)<\infty.\]
Since by
Lemma \ref{distance}(a) we have $s(k)\in X$, the relation
$\sum_{i=1}^m\left(f_i(s(k))-f_i(z^*)\right)\ge 0$ holds for all
$k$, thus implying that
\[\frac{1}{2}\sum_{k=0}^\infty \sum_{i=1}^m
\|\phi^i(k)\|^2<\infty,\]
\[\sum_{k=0}^\infty\a_k\left(f(s(k))-f(z^*)\right)<\infty.\]
In view of the former of the preceding two relations, we have
\[\lim_{k\to\infty} \phi^i(k)=0
\qquad\hbox{for all $i$},\] while from the latter, since
$\sum_{k}\a_k=\infty$ and $f(s(k))-f^*\ge0$ [because $s(k)\in X$ for
all $k$], we obtain \begin{equation}\liminf_{k\to\infty}
f(s(k))=f^*.\label{con-seq}\end{equation}

Since $\phi^i(k)\to0$ for all
$i$ and $\a_k\to 0$ [in view of $\sum_{k}\a_k^2<\infty$],
from Eq.\ (\ref{projdis}) it follows that
\[\lim_{k\to\infty}\|v^i(k)-s(k)\|=0\qquad\hbox{for all $i$}.\]
Finally, since $x^i(k+1)=  v^i(k) - \a_k d_i(k) + \phi^i(k)$ [see
(\ref{error_form})], in view of $\a_k\to0$, $\|d_i(k)\|\le L$, and
$\phi^i(k)\to0$, we see that $\lim_{k\to\infty}\|x^i(k+1) -
v^i(k)\|=0$ for all $i$. This and the preceding relation yield
\[\lim_{k\to\infty}\|x^i(k+1)-s(k)\|=0\qquad\hbox{for all $i$}.\]

We now show that the sequences $\{x^i(k)\},i=1,\ldots, m,$ converge
to the same limit point, which lies in the optimal solution set $
X^*$. By taking limsup as $N\to\infty$ in relation (\ref{eqn:ok})
and then liminf as $K\to\infty$, (while dropping the nonnegative
terms on the right hand side there), since $\sum_{k}\a_k^2<\infty$,
we obtain for any $z^*\in X^*,$
\[\limsup_{N\to\infty} \sum_{i=1}^m\|x^i(N+1)-z^*\|^2
\le \liminf_{K\to\infty} \sum_{i=1}^m\|x^i(K)-z^*\|^2,\] implying
that the scalar sequence $\left\{\sum_{i=1}^m\|x^i(k)-z^*\|\right\}$
is convergent for every $z^*\in X^*$. Since $\|x^i(k+1) -
s(k)\|\to0$ for all $i$, it follows that the scalar sequence
$\{\|s(k)-z^*\|\}$ is also convergent for every $z^*\in X^*$. In
view of $\liminf_{k\to\infty} f(s_k)=f^*$ [cf.\ Eq.\
(\ref{con-seq})], it follows that one of the limit points of
$\{s_k\}$ must belong to $X^*$; denote this limit by $x^*$. Since
$\{\|s(k)-z^*\|\}$ is convergent for $z^*=x^*$, it follows that
$\lim_{k\to\infty} s(k)=x^*.$ This and $\|x^i(k+1) - s(k)\|\to0$ for
all $i$ imply that each of the sequences $\{x^i(k)\}$ converges to a
vector $x^*$, with $x^*\in X^*$.
\end{proof}

\section{Conclusions}\label{conclusions}

We studied constrained consensus and optimization problems where
agent $i$'s estimate is constrained to lie in a closed convex set
$X_i$. For the constrained consensus problem, we presented a
distributed projected consensus algorithm and studied its
convergence properties. Under some assumptions on the agent weights
and the connectivity of the network, we proved that each of the
estimates converge to the same limit, which belongs to the
intersection of the constraint sets $X_i$. We also showed that the
convergence rate is geometric under an interior point assumption for
the case when agent weights are time-invariant and uniform. For the
constrained optimization problem, we presented a distributed
projected subgradient algorithm. We showed that with a stepsize
converging to zero fast enough, the estimates generated by the
subgradient algorithm converges to an optimal solution for the case
when all agent constraint sets are the same and when agent weights
are time-invariant and uniform.

The framework and algorithms studied in this paper motivate a number
of interesting research directions. One interesting future direction
is to extend the constrained optimization problem to include both
local and global constraints, i.e., constraints known by all the
agents. While global constraints can also be addressed using the
``primal projection" algorithms of this paper, an interesting
alternative would be to use ``primal-dual" subgradient algorithms,
in which dual variables (or prices) are used to ensure feasibility
of agent estimates with respect to global constraints. Such
algorithms have been studied in recent work \cite{saddle} for
general convex constrained optimization problems (without a
multi-agent network structure).

Moreover, in this paper, we presented convergence results for the
distributed subgradient algorithm for two cases: agents have
time-varying weights but the same constraint set; and agents have
time-invariant uniform weights and different constraint sets. When
agents have different constraint sets, the convergence analysis
relies on an error bound that relates the distances of the iterates
(generated with constant uniform weights) to each $X_i$ with the
distance of the iterates to the intersection set under an interior
point condition (cf.\ Lemma \ref{distance}). This error bound is
also used in establishing the geometric convergence rate of the
projected consensus algorithm with constant uniform weights. These
results can be extended using a similar analysis once an error bound
is established for the general case with time-varying weights. We
leave this for future work.

\newpage

\bibliographystyle{amsplain}
\bibliography{constrained_consensus}

\providecommand{\bysame}{\leavevmode\hbox to3em{\hrulefill}\thinspace}
\providecommand{\MR}{\relax\ifhmode\unskip\space\fi MR }
% \MRhref is called by the amsart/book/proc definition of \MR.
\providecommand{\MRhref}[2]{%
  \href{http://www.ams.org/mathscinet-getitem?mr=#1}{#2}
}
\providecommand{\href}[2]{#2}
\begin{thebibliography}{10}

\bibitem{aronszajn}
N.~Aronszajn, \emph{Theory of reproducing kernels}, Transactions of the
  American Mathematical Society \textbf{68} (1950), no.~3, 337--404.

\bibitem{book}
D.P. Bertsekas, A.~Nedi\'c, and A.E. Ozdaglar, \emph{Convex analysis and
  optimization}, Athena Scientific, Cambridge, Massachusetts, 2003.

\bibitem{distbook}
D.P. Bertsekas and J.N. Tsitsiklis, \emph{Parallel and distributed computation:
  Numerical methods}, Athena Scientific, Belmont, MA, 1989.

\bibitem{Blatt07}
D.~Blatt, A.~O. Hero, and H.~Gauchman, \emph{A convergent incremental gradient
  method with constant stepsize}, SIAM Journal of Optimization \textbf{18}
  (2007), no.~1, 29--51.

\bibitem{multiagent}
V.D. Blondel, J.M. Hendrickx, A.~Olshevsky, and J.N. Tsitsiklis,
  \emph{Convergence in multiagent coordination, consensus, and flocking},
  Proceedings of IEEE CDC, 2005.

\bibitem{boyd}
S.~Boyd, A.~Ghosh, B.~Prabhakar, and D.~Shah, \emph{Gossip algorithms: Design,
  analysis, and applications}, Proceedings of IEEE INFOCOM, 2005.

\bibitem{spielman}
M.~Cao, D.A. Spielman, and A.S. Morse, \emph{A lower bound on convergence of a
  distributed network consensus algorithm}, Proceedings of IEEE CDC, 2005.

\bibitem{deutsch}
F.~Deutsch, \emph{Rate of convergence of the method of alternating
  projections}, Parametric Optimization and Approximation (B.~Brosowski and
  F.~Deutsch, eds.), vol.~76, Birkhäuser, Basel, 1983, pp.~96--107.

\bibitem{deuhun}
F.~Deutsch and H.~Hundal, \emph{The rate of convergence for the cyclic
  projections algorithm i: Angles between convex sets}, Journal of
  Approximation Theory \textbf{142} (2006), 36--55.

\bibitem{facc_pang}
F.~Facchinei and J-S. Pang, \emph{Finite-dimensional variational inequalities
  and complementarity probems}, Springer-Verlag New York, Vol.\ 1, 2003.

\bibitem{gubinpolyak}
L.G. Gubin, B.T. Polyak, and E.V. Raik, \emph{The method of projections for
  finding the common point of convex sets}, {U.S.S.R} Computational Mathematics
  and Mathematical Physics \textbf{7} (1967), no.~6, 1211--1228.

\bibitem{ali}
A.~Jadbabaie, J.~Lin, and S.~Morse, \emph{Coordination of groups of mobile
  autonomous agents using nearest neighbor rules}, IEEE Transactions on
  Automatic Control \textbf{48} (2003), no.~6, 988--1001.

\bibitem{Johansson07}
B.~Johansson, M.~Rabi, and M.~Johansson, \emph{A simple peer-to-peer algorithm
  for distributed optimization in sensor networks}, Proceedings of the 46th
  IEEE Conference on Decision and Control, 2007, pp.~4705--4710.

\bibitem{JasonECC}
J.R. Marden, G.~Arslan, and J.S. Shamma, \emph{Connections between cooperative
  control and potential games illustrated on the consensus problem},
  Proceedings of the 2007 European Control Conference, 2007.

\bibitem{Nedic01}
A.~Nedi\'{c} and D.~P. Bertsekas, \emph{Incremental subgradient method for
  nondifferentiable optimization}, SIAM Journal of Optimization \textbf{12}
  (2001), 109--138.

\bibitem{alex_quant_pap}
A.~Nedi\'c, A.~Olshevsky, A.~Ozdaglar, and J.~N. Tsitsiklis, \emph{On
  distributed averaging algorithms and quantization effects}, LIDS Technical
  Report 2778, available at http://arxiv.org/abs/0803.1202, 2007.

\bibitem{ratesubgrad}
A.~Nedi\'c and A.~Ozdaglar, \emph{On the rate of convergence of distributed
  subradient methods for multi-agent optimization}, Proceedings of IEEE CDC,
  2007.

\bibitem{ratejournal}
\bysame, \emph{Distributed subradient methods for multi-agent optimization},
  IEEE Transactions on Automatic Control, forthcoming, 2008.

\bibitem{saddle}
\bysame, \emph{Subgradient methods for saddle-point problems}, Journal of
  Optimization Theory and Applications, forthcoming, 2008.

\bibitem{vonNeumann}
J.~Von Neumann, \emph{Functional operators}, Princeton University Press,
  Princeton, 1950.

\bibitem{reza}
R.~Olfati-Saber and R.M. Murray, \emph{Consensus problems in networks of agents
  with switching topology and time-delays}, IEEE Transactions on Automatic
  Control \textbf{49} (2004), no.~9, 1520--1533.

\bibitem{alexCDC}
A.~Olshevsky and J.N. Tsitsiklis, \emph{Convergence rates in distributed
  consensus averaging}, Proceedings of IEEE CDC, 2006.

\bibitem{alexlong}
\bysame, \emph{Convergence speed in distributed consensus and averaging}, SIAM
  Journal on Control and Optimization, forthcoming, 2008.

\bibitem{Sundhar08a}
S.~Sundhar Ram, A.~Nedi\'c, and V.~V. Veeravalli, \emph{Incremental stochastic
  sub-gradient algorithms for convex optimization}, Available at
  http://arxiv.org/abs/0806.1092, 2008.

\bibitem{Rock70}
R.~T. Rockafellar, \emph{Convex analysis}, Princeton University Press, 1970.

\bibitem{johnthes}
J.N. Tsitsiklis, \emph{Problems in decentralized decision making and
  computation}, Ph.D. thesis, {M}assachusetts {I}nstitute of {T}echnology,
  1984.

\bibitem{distasyn}
J.N. Tsitsiklis, D.P. Bertsekas, and M.~Athans, \emph{Distributed asynchronous
  deterministic and stochastic gradient optimization algorithms}, IEEE
  Transactions on Automatic Control \textbf{31} (1986), no.~9, 803--812.

\bibitem{vicsek}
T.~Vicsek, A.~Czirok, E.~Ben-Jacob, I.~Cohen, and Schochet O., \emph{Novel type
  of phase transitions in a system of self-driven particles}, Physical Review
  Letters \textbf{75} (1995), no.~6, 1226--1229.

\bibitem{wolf}
J.~Wolfowitz, \emph{Products of indecomposable, aperiodic, stochastic
  matrices}, Proceedings of the American Mathematical Society \textbf{14}
  (1963), no.~4, 733--737.

\end{thebibliography}

\end{document}